\documentclass[12pt]{article}
\usepackage{times}

\usepackage[latin1]{inputenc}

\usepackage{amsmath}
\usepackage{amssymb}
\usepackage{mathrsfs}
\usepackage{xspace}
\usepackage{units}              % for "nicefrac"

\usepackage[text={5.5in, 9in}, centering]{geometry}

\usepackage[ps,tips,matrix,arrow]{xy}
\SelectTips{times}{}

\ifx\pdftexversion\undefined
  \usepackage[dvips]{graphicx}
\else
  \usepackage[pdftex]{graphicx}
\fi

\usepackage[amsmath,thmmarks]{ntheorem}

\newcommand{\inlinexy}[1]{{\def\objectstyle{\scriptstyle}\def\labelstyle{\scriptstyle}{\xymatrix@-1.2pc@1{#1}}}}

\newcommand{\defsep}{\nopagebreak\medskip}

\newcommand{\comp}{\circ}       % Function composition
\newcommand{\R}{{\mathbb R}}        % The reals
\newcommand{\N}{{\mathbb N}}        % The naturals

\newcommand{\set}[1]{\left\{#1\right\}}

\newcommand{\tuple}[1]{\left(\smash{#1}\right)}
\newcommand{\gens}[1]{\left<\smash{#1}\right>}

\newcommand{\bdry}{\partial}
\newcommand{\inv}{^{-1}}
\newcommand{\tensor}{\otimes}
\newcommand{\bigtensor}{\bigotimes}

\newcommand{\map}{\operatorname{Map}}

\newcommand{\cross}{\times}
\newcommand{\homeo}{\cong}
\newcommand{\isom}{\cong}

\newcommand{\restr}[2]{{#1}|_{#2}}
\newcommand{\inclto}{\hookrightarrow}

\newcommand{\projto}{\twoheadrightarrow}

\newcommand{\from}{\leftarrow}
\newcommand{\xfrom}{\xleftarrow}
\newcommand{\xto}{\xrightarrow}

\newcommand{\xinclfrom}[1]{\stackrel{#1}\hookleftarrow}

\newcommand{\disjcup}{\sqcup}
\newcommand{\bigdisjcup}{\bigsqcup}
\renewcommand{\emptyset}{\varnothing}
\newcommand{\ob}{\operatorname{Ob}}

\newcommand{\blank}{\text{\bf ---}\hspace{1pt}}
\renewcommand{\hom}{\operatorname{Hom}}
\newcommand{\sym}{\operatorname{Sym}} %symmetric group on a set
\newcommand{\aut}{\operatorname{Aut}} %automorphisms

\newcommand{\id}{\mathrm{id}}

\newcommand{\downto}{\searrow}

\newcommand{\ev}{\mathrm{ev}} % Evaluation map 
\newcommand{\pt}{\ast}             % point
\newcommand{\half}{{\nicefrac12}}
\newcommand{\B}{{\mathcal B}}

\newcommand{\M}{{\mathcal M}}
\newcommand{\G}{\Gamma}
\newcommand{\tG}{\tilde{\Gamma}}

\newcommand{\V}{{\mathcal V}}

\newcommand{\bdryS}{\bdry_{S^1}}
\newcommand{\bdryI}{\bdry_I}
\newcommand{\bdryW}{\bdry_{\mathrm{w}}}
\newcommand{\bdryF}{\bdry_{\mathrm{f}}}
\newcommand{\bdrySt}{\bdry_{\mathrm{s}}}
\newcommand{\iotaF}{\iota_{\mathrm{f}}}
\newcommand{\iotaSt}{\iota_{\mathrm{s}}}

\newcommand{\sop}{\mu}
\renewcommand*{\P}[2]{P{#1}\tuple{#2}}

\newcommand{\glue}{\text{\tt \#\,}}
\newcommand{\Gg}{\G_1\glue\G_2}
\newcommand{\Fat}{{\mathscr F}\!\!\textit{at}}
\newcommand{\FatP}{{\mathscr F}\!\!\textit{at}^\star}

\newcommand{\aFP}{{\mathscr F}^+}

\newcommand{\Graph}{\mathscr{G}\!\textit{raph}}
\renewcommand{\phi}{\varphi}
\newcommand{\mM}{{\mathbb M}}   
   
\newcommand{\fV}{{\mathscr V}}   
\newcommand{\catC}{{\mathscr C}}        % Category C
\newcommand{\pCob}{\mathrm{Cob}}
\newcommand{\pCobpb}{\mathrm{Cob}^+}
\newcommand{\op}{\mathrm{op}}   % Opposite

\newcommand{\bmap}[2]{[\![#1]\!]_{#2}}
\newcommand{\bmapm}[1]{\bmap{#1}{\M}}
\newcommand{\bma}[1]{[\![#1]\!]}
\newcommand{\bsigma}{\hat{\sigma}}

\newcommand{\tqft}{{\scshape tqft}\xspace}
\newcommand{\tqfts}{{\scshape tqft}s\xspace}
\newcommand{\prop}{{\scshape prop}\xspace}
\newcommand{\props}{{\scshape prop}s\xspace}
\newcommand{\bprop}{$\B$\nobreakdash-{\scshape prop}\xspace}

\newcommand{\btqft}{{$\B$-{\scshape tqft}}\xspace}
\newcommand{\btqfts}{$\B$-{\scshape tqft}s\xspace}

\newcommand{\bcobos}{$\B$-cobordisms\xspace}
\newcommand{\oc}{open-closed\xspace}
\newcommand{\Oc}{Open-closed\xspace}
\newcommand{\ocobo}{open-closed cobordism\xspace}
\newcommand{\ocobos}{open-closed cobordisms\xspace}
\newcommand{\osurf}{open-closed surface\xspace}
\newcommand{\osurfs}{open-closed surfaces\xspace}
\newcommand{\bfam}{$\B$-family\xspace}

\newcommand{\bgraphs}{$\B$-graphs\xspace}
\newcommand{\bgraph}{$\B$-graph\xspace}
\newcommand{\fgraph}{fat graph\xspace}
\newcommand{\fgraphs}{fat graphs\xspace}

\newcommand{\ofgraph}{\oc fat graph\xspace}

\newcommand{\ofgraphs}{\oc fat graphs\xspace}
\newcommand{\Ofgraphs}{\Oc fat graphs\xspace}

\newcommand{\Osurfs}{Open-closed surfaces\xspace}

\newcommand{\Fgraphs}{Fat graphs\xspace}

\newcommand{\horiented}{$h_*$\nobreakdash-oriented\xspace}

\newcommand{\graphdef}[3]{\begin{itemize}
  \item \emph{Vertices:} #1\nopagebreak
  \item \emph{Half-edges:} #2
  \item \emph{Incidence of half-edges:} #3
  \end{itemize}}

\newsavebox{\thevcbox} 
\newlength{\vcboxwidth}
\newcommand{\vcbox}[1]{
  \savebox{\thevcbox}{#1}
  \settowidth{\vcboxwidth}{\usebox{\thevcbox}}
  \begin{minipage}[c]{\vcboxwidth}\usebox{\thevcbox}\end{minipage}
}

\newcommand{\picbox}[1]{\vcbox{\includegraphics{#1}}}

\newcommand{\inclar}{\ar@{^(->}}
\newcommand{\surjar}{\ar@{->>}}
\newcommand{\injar}{\ar@{>->}}

\newenvironment{proofof}[1]{\begin{proofNoheader}[Proof of #1]}{\end{proofNoheader}}

\begin{document}
\title{Open-Closed String Topology via Fat Graphs} \author{Antonio
  Ramírez\thanks{The author was supported by Stanford University, the
    Pacific Institute of Mathematical Sciences at the University of
    British Columbia, and the Mathematical Sciences Research Institute
    throughout the preparation of this article. }}
    \maketitle

\begin{abstract}
  Given a smooth closed manifold $M$ with a family $\set{L_i}$ of
  closed submanifolds, we consider the free loop space $LM$ and
  the spaces $\P{M}{L_i,L_j}$ of open strings (paths $\gamma:[0,1]\to
  M$ with $\gamma(0)\in L_i,$ $\gamma(1)\in L_j$).  We construct
  string topology operations resulting in an \emph{open-closed} \tqft
  on the family $(h_*(LM),\set{h_*(\P{M}{L_i,L_j})}_{i,j\in\B})$ which
  extends the known string topology \tqft on $h_*(LM)$.  Here, $h_*$
  is a multiplicative generalized homology theory supporting
  orientations for $M$ and the $L_i$.  To construct the operations, we
  introduce the notion of \emph{\ofgraph,} generalizing \fgraphs to
  the open-closed setting.
\end{abstract}

\section{Introduction}
\label{sec:introduction}
The area of string topology began with a construction by M.~Chas and
D.~Sullivan~\cite{chas-sullivan:string-topology} of previously
undiscovered algebraic structure on the homology of the free loop
space of a closed oriented manifold $M.$ This is is the space $LM$ of
all continuous maps from the circle to $M.$ Among other results, Chas
and Sullivan found that the homology of $LM,$ with its grading
suitably shifted, is a graded commutative algebra, much like the
homology of an oriented manifold does by virtue of Poincaré duality.

This operation, the \emph{string loop product,} can be understood by
considering the space $\map(P,M)$ of maps from a pair-of-pants surface
$P$ to $M,$ where $P$ is regarded as a cobordism from a disjoint union
of two circles to a single circle. $P$ is a model for the basic
interaction in string theory, in which two strings merge to form a
single string.

There is a diagram
\begin{equation}\label{eq:1}%\tag{$\star$}
LM\cross LM \xfrom{i} \map(P,M) \xto{j} LM,
\end{equation}
where~$i$ and~$j$ are restriction maps to the incoming and outgoing
boundary of~$P.$ The product is then constructed as the composition
\begin{displaymath}
  H_*(LM)\tensor H_*(LM) \xto{i^!}
H_*(\map(P,M)) \xto{j_*} H_*(LM),
\end{displaymath}
where $i^!$ is an umkehr map.

Since the spaces involved are infinite-dimensional, the existence of
$i^!$ is not immediate. However, it can be constructed by observing
that $P$ is homotopy equivalent to a ``figure-eight'' space
$\G=S^1\vee S^1,$ so that we may replace (\ref{eq:1}) by
\begin{equation}\label{eq:2}%\tag{$\star\star$}
LM\cross LM \xfrom{i} \map(\G,M) \xto{j} LM,
\end{equation}
where the maps are finite codimension embeddings in an appropriate
sense. Using essentially diagram~(\ref{eq:2}), Chas and Sullivan
constructed the operation using transversality of smooth chains in
$LM.$ There have been other constructions since then, notably
including the homotopy-theoretic approach by R.~Cohen and
J.~Jones~\cite{cohen:htpy-theoretic-realization} in which $i^!$ is
defined using a Pontrjagin-Thom collapse.

\medskip R.~Cohen and
V.~Godin~\cite{cohen-godin:polarized-view-string-topo} showed later
that the pair-of-pants $P$ may be replaced by any oriented connected
cobordism $\Sigma$ between one-manifolds, provided that it has
nonempty outgoing boundary. The result is a family of operations
\begin{displaymath}\mu_\Sigma:~H_*(LM)^{\tensor p}\to H_*(LM)^{\tensor q}\end{displaymath} 
compatible with gluing of cobordisms. This is a form of
\emph{topological quantum field theory} (\tqft). The construction
exploits the fact that any surface with boundary may be represented as
a \emph{\fgraph,} which is a finite graph $\G$ endowed with extra data
that determines a surface with boundary having $\G$ as a
deformation retract.

\medskip The appearance of \fgraphs reflects a result, due in its
various forms to Harer \cite{harer:virtual-coho-dim-mcg},
Penner~\cite{penner:decorated-teichmuller-space}, and
Strebel~\cite{strebel:quadratic-differentials}, that moduli spaces of
punctured Riemann surfaces are homotopy equivalent to spaces of metric
\fgraphs.

\subsection{Open-closed string topology}
\emph{Open-closed string topology} generalizes string topology by
allowing the strings to be ``open,'' that is, paths in a manifold $M$
which need not be loops.  The endpoints of open strings are
constrained to lie in certain distinguished submanifolds $L_b\subseteq
M$. The basic spaces of open strings considered are then of the form
$\P{M}{L_1,L_2},$ standing for the space of paths $\gamma:[0,1]\to M$
such that $\gamma(0)\in L_1$ and $\gamma(1)\in L_2.$ These
submanifolds (or typically objects with more structure) are called
\emph{$D$-branes} in string theory, and we will shorten it to
``branes.''

The idea of open string topology was introduced by Sullivan in
\cite{sullivan:openclosed}, using transversality of smooth chains. The
prototype construction, in homotopy-theoretic terms, is as follows.
Consider the diagram
\begin{displaymath}\xymatrix{
  {\P{M}{L_1,L_2}\cross \P{M}{L_2,L_3}\ar[d]_{\ev_1\cross\ev_0}}   & {\P{M}{L_1,L_2,L_3}\ar[r]^-{j}\ar[l]_-{i}\ar[d]_{\ev_{\half}}} & {\P{M}{L_1,L_3}} \\
  {L_2\!\cross\! L_2}                          & {L_2\ar[l]_{\Delta}.} &
}\end{displaymath}
Here, we define $\P{M}{L_1,L_2,L_3}=\set{\gamma\in \P{M}{L_1,L_3}:
  \gamma(\half)\in L_2}$. The map $i$ takes a path in
$\P{M}{L_1,L_2,L_3}$ and splits it at the middle into two paths, $j$
is the inclusion, and the vertical maps are evaluation maps.  We may
define a \emph{composition} operation
\begin{displaymath}H_*(\P{M}{L_1,L_2})\tensor H_*(\P{M}{L_2,L_3})\to H_*(\P{M}{L_1,L_3})\end{displaymath}
as~$j_*\comp i^!$. When $M$ and $L_2$ are oriented, the umkehr map
$i^!$ exists by the Pontrjagin-Thom collapse, because~$i$ is a finite
codimension inclusion in a suitable sense. See
Section~\ref{sec:thom-collapse} for a general statement.

\bigskip Our aim is to extend the string topology \tqft of
\cite{cohen-godin:polarized-view-string-topo} to the \oc
setting.  Here the one-manifolds considered may have endpoints, which
carry boundary conditions. Thus each endpoint is labeled by an element
of a set $\B$ indexing a family $\set{L_b}_{b\in\B}$ of branes.
Cobordisms between closed one-manifolds are replaced accordingly by
the natural cobordisms between one-manifolds with labeled boundary.
We call the resulting structures \emph{\oc} \tqfts, or \btqfts to make
explicit the dependence on $\B$; see Definition~\ref{def:btqft}. The
main result is as follows.
\begin{theoremNoheader}[Theorem~\ref{thm:main1}]
  Let $M$ be a closed smooth manifold, and let $\set{L_b}_{b\in\B}$ be
  a family of smooth closed submanifolds. Suppose that $h_*$ is a
  multiplicative generalized homology theory whose coefficient ring
  $h_*(\pt)$ is a graded field, and suppose that $M$ and the $L_b$ are
  oriented with respect to $h_*$. Then, the family
  \begin{displaymath}
    \tuple{h_*LM,\set{h_*\P{M}{L_a,L_b}}_{a,b\in\B}}
  \end{displaymath}
  supports a positive-boundary \btqft structure over the coefficient
  ring.  This extends the known string topology \tqft on~$h_*(LM)$.
\end{theoremNoheader}
The qualifier ``positive boundary'' means that \btqft operations only
exist for those cobordisms having nonempty outgoing boundary on each
connected component. This restriction is also present in the closed
case.

\bigskip The main tools for constructing the operations will be
{\emph\ofgraphs}~(Definition~\ref{def:ofgraph}). \Ofgraphs are a
generalization of \fgraphs. They have a well-defined notion of
``string boundary,'' which is a graph isomorphic to a one-manifold
with $\B$-labeled boundary. They also have a ``fattening'' operation,
which yields an associated ``\osurf.'' These two properties generalize
the corresponding properties of \fgraphs.

\bigskip The content of this paper is essentially part of the author's
PhD thesis \cite{ramirez:thesis}. We have recently learned of work by
E.~Harrelson~(\cite{harrelson:homology-openclosed},
\cite{harrelson:thesis}) that overlaps with the present paper. We hope
that this will be the subject of future collaboration.

\medskip We are very thankful to Professor Ralph L. Cohen for his
great patience and generous guidance throughout this project. 

% --------------------------------------------------------------------
\section{\Osurfs and cobordisms}
\label{sec:osurfs}
% --------------------------------------------------------------------

Assume given an arbitrary set $\B$ of ``formal branes,'' which we will
eventually use to index the actual branes.

\begin{definition}\label{def:labeled-1-manifold} \label{def:open-closed-surface}
  A \emph{$\B$-labeled one-manifold} $C$ is an oriented one-manifold
  with boundary together with a function $\beta:\bdry C\to\B$ (the
  $\B$-labeling). An \emph{isomorphism} between two $\B$-labeled
  one-manifolds is a diffeomorphism that preserves the orientation and
  the $\B$-labeling. Let $C^*$ stand for $C$ with the opposite
  orientation. Given $a, b\in\B$, let $I_{a,b}$ be a copy of the unit
  interval, oriented in the direction from $0$ to $1$, with
  $\beta(0)=a$, $\beta(1)= b$; note that $I_{a,b}^*$ is isomorphic to
  $I_{b,a}$.

  \defsep An \emph{\osurf} (with brane labels drawn from $\B$) is a smooth
  oriented surface with boundary $S$, together with distinguished
  embedded one-dimensional submanifolds with boundary $\bdrySt S,
  \bdryF S\subseteq \bdry S$ (the \emph{string boundary} and
  \emph{free boundary,} respectively) and a locally constant function
  $\bdryF S \to\B$ (the \emph{brane labeling}) such that:
  \begin{enumerate}
  \item $\bdry S = \bdrySt S \cup \bdryF S$, and,
  \item $\bdry(\bdrySt S) = \bdry(\bdryF S) = \bdrySt S \cap \bdryF
    S$.
  \end{enumerate}
  The restriction of $\beta$ to $\bdry(\bdrySt S)$ makes $\bdrySt S$ a
  $\B$-labeled one-manifold.
\end{definition}

\begin{figure}
  \centering
  \fbox{\begin{tabular}{cc}
    \picbox{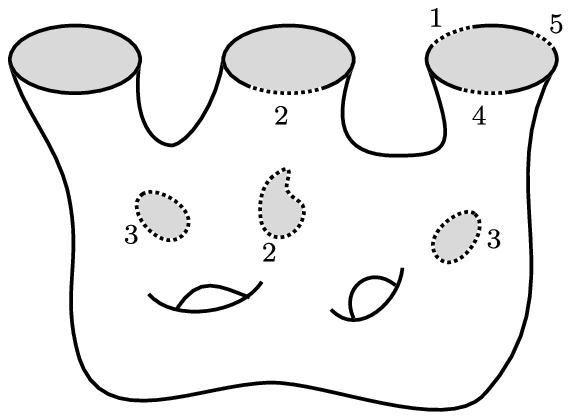}\!\!\!\!\!\!\!\!\! &   \!\!\!\picbox{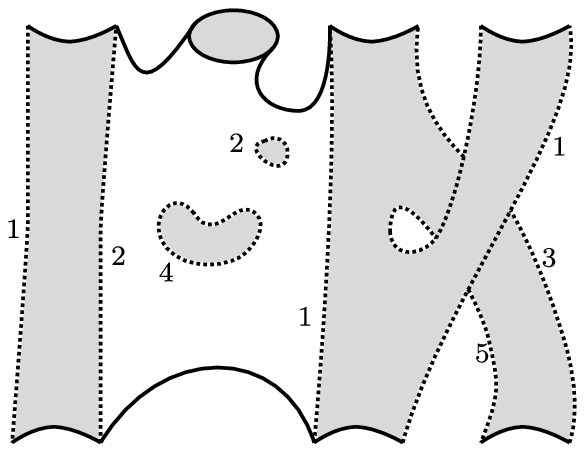} \\
    (a) & (b)
  \end{tabular}}  \caption{ With $\B=\set{1,2,3,\ldots}$: (a) An \osurf. (b) An \ocobo from $I_{1,2}\disjcup I_{1,1}\disjcup I_{5,3}$ to $I_{1,2}\disjcup S^1\disjcup I_{1,3}\disjcup I_{5,1}$.}
  \label{fig:osurf}
\end{figure}

\begin{example}
  Figure~\ref{fig:osurf}(a) shows an \osurf with genus two, five
  string boundary components (only one of which is closed), and seven
  free boundary components, drawn dotted, with brane labeling in
  $\B=\set{1,2,\ldots}$.
\end{example}

\begin{definition}\label{def:ocobo}
  Given $\B$-labeled one-manifolds $C_-$, $C_+$, an \emph{\ocobo}
  from $C_-$ to $C_+$ is a $\B$-surface $S$ together with a
  decomposition $\bdrySt S=\bdry^-S\disjcup\bdry^+S$ and
  orientation-preserving isomorphisms $\iota_-:C_-\to\bdry^- S$,
  $\iota_+:C_+^*\to\bdry^+ S$ of $\B$-labeled one-manifolds. We call
  $\bdry^-S$, $\bdry^+S$ the \emph{incoming} and \emph{outgoing}
  boundary, respectively.
\end{definition}

\begin{example}
  Figure~\ref{fig:osurf}(b) shows an \ocobo $S$, where
  the incoming boundary is at the bottom and the outgoing boundary at
  the top.  Notice that a single topological boundary component may
  contain both incoming and outgoing string boundary subintervals.
\end{example}

% --------------------------------------------------------------------
\section{Open-closed TQFT}\label{sec:btqft}
% --------------------------------------------------------------------
Recall that an ordinary \tqft consists of a vector space $V$ together
with a homomorphism $\mu_\Sigma: V^{\tensor p}\to V^{\tensor q}$ for
every oriented cobordism $\Sigma$ from $\bigdisjcup_p S^1 $ to
$\bigdisjcup_q S^1$. The maps $\mu_\Sigma$ are required to be
diffeomorphism-invariant and to satisfy natural compatibility
conditions with respect to disjoint union and composition of
cobordisms. This can be described succinctly using the language of
\props: a \tqft is an algebra over $\pCob$, which is a \prop of
oriented cobordisms between circles (see, for example, Voronov
~\cite{voronov:univ-algebra-notes}). We will introduce an analogous
notion to describe \tqfts involving open and closed strings.

\begin{definition}
  Define $\mM_\B$ to be the free abelian monoid on the symbols $S^1$
  and $I_{a,b}$ for $a,b\in\B$. Define a \emph{\bprop} to be a
  symmetric strict monoidal category having $\mM_\B$ as its monoid of
  objects. 
\end{definition}

This can be elaborated as follows. The symmetric monoidal
axioms~\cite{maclane:categories} imply that each object
$$x=nS^1+\sum_{a,b\in B} m_{a,b}I_{a,b}\in\mM_\B$$ 
of a \bprop $\catC$ carries an action $\Sigma_x\to\aut_{\catC}(x)$,
where the group
$$\Sigma_x = \Sigma_n\cross\prod_{a,b\in B}\Sigma_{m_{a,b}}$$ is a
permutation group associated to $x$. It follows that each set
$\catC(x,y)$ has an action by $\Sigma_x$ on the right and a commuting
action by $\Sigma_y$ on the left.  Moreover, for any
$x,y,z\in\ob(\catC)$, the composition map
$\comp:\catC(x,y)\cross\catC(y,z)\to\catC(x,z)$ is equivariant with
respect to the $\Sigma_x$ and $\Sigma_z$ actions, and it satisfies
$p\comp(\sigma q)=(p\sigma)\comp q$ for every $p\in\catC(y,z)$,
$q\in\catC(x,y)$, $\sigma\in\Sigma_y$.

% Notice the abuse of notation under which $I_{a,b}$ may stand for a
% formal symbol in $\mM_\B$ or for a $\B$-labeled unit interval;
% likewise for $S^1$.

\begin{definition}
  An \emph{algebra} over a \bprop $\catC$ is defined to be a monoidal
  functor from $\catC$ to the symmetric monoidal category of
  $R$-modules for some ring $R$.
\end{definition}
In more detail, an algebra over $\catC$ is specified by giving:
\begin{enumerate}
\item a \emph{\bfam} $\fV=(V,\set{V_{a,b}}_{a,b}\in\B)$ of $R$-modules, and,
\item for every $x,y\in\mM_\B$, a map $\catC(x,y)\to
  \hom(\fV(x),\fV(y))$, where we define
  \begin{displaymath}
    \fV(nS^1+\sum_{a,b\in\B}m_{a,b}I_{a,b}):=V^{\tensor
      n}\tensor\bigtensor_{a,b\in B} V_{a,b}^{\tensor m_{a,b}}.
  \end{displaymath}
\end{enumerate}
These maps are required to satisfy the necessary functoriality and
equivariance conditions. 
% The sets $\hom(\fV(x),\fV(y))$ are the
% morphism sets of the \emph{endomorphism~\bprop} $\End(\fV)$ of the
% \bfam. In these terms, an algebra over $\catC$ may be described as a
% morphism of \bprops from $\catC$ to $\End(\fV)$.

\begin{definition}\label{def:pcobB}
  Given any $x=nS^1+\sum_{a,b\in\B}m_{a,b}I_{a,b}\in\mM_\B$, let $|x|$
  be the $\B$-labeled one-manifold $(S^1)^{\disjcup
    n}\disjcup\bigdisjcup_{a,b\in \B} I_{a,b}^{\disjcup m_{a,b}}$ (in
  particular, $|0|=\emptyset$).

  \defsep%
  Define the \bprop $\pCob_\B$ (the \emph{cobordism} \bprop) by
  letting a morphism from $x$ to $y$ be an equivalence class of
  triples $(S,L_-,L_+)$, where:
  \begin{enumerate}
  \item $S$ is a $\B$-cobordism with $\bdry^-S\isom |x|$ and
    $\bdry^+S\isom |y|^*$, together with a choice of parametrization
    for each string boundary component by either $I$ or $S^1$; this
    parametrization is orientation-preserving on the incoming components
    and orientation-reversing on the outgoing components.
  \item $L_-$ is a function $\pi_0(\bdry^-S)\to\N$ which restricts to
    a bijection with $\set{1,\ldots,m_{a,b}}$ on incoming components of
    type $I_{a,b}$ and to a bijection with $\set{1,\ldots,n}$ on incoming
    components of type $S^1$. Similarly for $L_+:\pi_0\bdry^+S\to\N$.
  \item Two such triples are considered equivalent if they are related
    by an isomorphism of $\B$-cobordisms which preserves the boundary
    parametrizations and the orderings $L_\pm$ of boundary components.
  \end{enumerate}

  The composition of $(S,L_-,L_+)\in\pCob_\B(x,y)$ and
  $(S',L'_-,L'_+)\in\pCob_\B(y,z)$ is given by
  $(S\cup_{|y|}S',L_-,L'_+)$, where $S\cup_{|y|}S'$ is the
  $\B$-cobordism that results from identifying each outgoing component
  $c$ of $\bdry^+S$ to the unique incoming component
  $c'\subset\bdry^-S'$ of the same type such that $L_+(c)=L_-(c')$.
  The identification is with respect to the boundary parametrizations
  which are part of the data.  The monoidal structure is given by
  letting $(S,L_-,L_+)\tensor(S',L'_-,L'_+):=(S\disjcup S',
  \tilde{L}_-,\tilde{L}_+)$, where each labeling $\tilde{L}_\pm$ is
  given by ordering the boundary components of $S'$ after those of $S$
  in each $\B$-labeling type. The group $\Sigma_y\cross\Sigma_x^\op$
  then acts on $(S,L_-,L_+)\in\pCob_\B(x,y)$ by letting $\Sigma_x$
  permute the labeling $L_-$ and letting $\Sigma_y$ permute $L_+$.
\end{definition}

\begin{definition} \label{def:btqft} An \emph{\oc topological
    quantum field theory with branes~$\B$} (which we will abbreviate
  \btqft) is an algebra over the \bprop $\pCob_\B$.
\end{definition}
This definition restricts to the usual definition of a \tqft when
$\B=\emptyset$.

\medskip The string topology operations do not yield a whole \btqft,
since there are no operations associated to \bcobos to the empty
one-manifold.  The appropriate variant is as follows.
\begin{definition} \label{def:pb-btqft} Let $\pCobpb_\B$ be the
  subcategory of $\pCob_\B$ consisting of \bcobos in which every
  connected component has nonempty outgoing boundary.  This category
  inherits a \bprop structure from $\pCob_\B$. A \emph{positive
    boundary} \btqft is an algebra over $\pCobpb_\B$.
\end{definition}

\begin{remark}
  We regard this definition chiefly as an ad-hoc device, and we do not
  claim that it is the most adequate definition of an open-closed
  field theory. For other treatments, we refer the reader to A. D.
  Lauda and H. Pfeiffer \cite{lauda-pfeiffer:openclosed} and to K.
  Costello \cite{costello:cy-cats}.
\end{remark}

% --------------------------------------------------------------------
\section{Graphs}
% --------------------------------------------------------------------

\Fgraphs have been used extensively to study moduli spaces of
punctured and bordered Riemann surfaces.  A \fgraph consists of a
finite graph together with a cyclic ordering of the edges (more
precisely, half-edges) incident at every given vertex. This extra
structure specifies a canonical ``fattening'' of the graph having the
form of an oriented surface with boundary with the graph as a
deformation retract, in which the punctures correspond to certain
cycles of oriented edges.

We will define the notion of \emph{\ofgraph,} which is a \fgraph with
extra structure, which will induce an \osurf structure on the
fattening. To fix notation, we recall some basic definitions.

\begin{definition}\label{def:graph}
  A \emph{graph} $\G=\tuple{V(\G),H(\G),s,r}$ consists of
  \begin{itemize}
  \item a set $V(\G)$ of \emph{vertices},
    
  \item a set $H(\G)$ of \emph{half-edges},
    
  \item a ``source'' map $s:H(\G)\to V(\G)$ taking a half-edge to the
    vertex that it attaches to, and,
    
  \item a ``reversal'' involution $r$ of $H(\G)$, having no fixed
    points, understood to take a half-edge to the opposite half-edge
    of their common edge.
  \end{itemize}
  The \emph{edges} $E(\G)$ are the orbits of $r$. We introduce the
  notation $H(v):=s\inv(v)$ for the set of half-edges incident with
  a given vertex. The \emph{degree} (or \emph{valence}) of a vertex
  $v\in V(\G)$ is defined to be the number $\#H(v)$.  

  \defsep%
  Define $t$ as the composition $s\comp r$, which is the
  ``target'' map taking a half-edge to its destination vertex (the
  source of its reversal).
  
  \defsep%
  Say that a graph is \emph{discrete} if it has no edges; any set (in
  particular, $\B$) can then be regarded as a discrete graph, and we
  will do so without mention.

  \defsep%
  Given two graphs $\G_1, \G_2$, define a morphism $\phi:\G_1\to\G_2$
  of graphs to be a pair $\tuple{\phi_V,\phi_H}$ of functions
  $\phi_V:V(\G_1)\to V(\G_2)$, $\phi_H: H(\G_1)\to H(\G_2)\disjcup
  V(\G_2)$, such that
  \begin{enumerate}
  \item $\phi_V(s(e))=s(\phi_H(e))$ for all $e\in H(\G_1)$, and,
  \item $\phi_H(r(e))=r(\phi_H(e))$ for all $e\in H(\G_1)$,
  \end{enumerate}
  where we extend the structure maps $s$ and $r$ to $V(\G_2)$ as the
  identity. When unambiguous, we will refer to both maps $\phi_V,
  \phi_H$ simply by $\phi$.  Define $\Graph$ to be the category of
  finite graphs, with these morphisms.

  \defsep Define a \emph{subgraph} of a graph
  $\G$ to be a graph $\G'$ with $H(\G')\subseteq H(\G)$,
  $V(\G')\subseteq V(\G)$ and for which the edge reversal and source
  maps are given by restriction from those of $\G$; we write
  $\G'\subseteq \G$.
\end{definition}
All our graphs (possibly excluding $\B$) will be finite.

\begin{remark}
  A half-edge $e\in H(\G)$ can be equally regarded as an oriented
  edge, oriented (for definiteness) in the direction that points away
  from its source vertex $s(e)$. We will use this point of view when
  convenient.
\end{remark}

Given $\G'\subseteq\G$, we would like a complement graph
$\G\setminus\G'$. This is obtained naively by removing from $\G$ all
vertices and edges belonging to $\G'$. After this, though, every $e\in
H(\G)\setminus H(\G')$ attached to a vertex $s(e)\in V(\G')$ ends up
with no source vertex.  We repair the result by formally introducing a
new vertex for each such $e$. Precisely:
\begin{definition} \label{def:subgraph-complement} Given
  $\G'\subseteq\G$, define the \emph{complement} graph
  $\G\setminus\G'$ as follows: \graphdef {%Vertices
    $V(\G\setminus\G') := (V(\G)\setminus V(\G'))\disjcup\delta(\G\setminus \G)',$
    where $\delta(\G\setminus\G'):=s_\G\inv(V(\G'))\setminus H(\G')$ is the set of
    half-edges attached to a vertex of $\G'$ but not lying in $\G'$. }
  {%Half-edges
    $H(\G\setminus\G') := H(\G)\setminus H(\G')$, with edge-reversal
    involution given by restriction from that of $\G$.  }{%Incidence
    Define $s_{\G\setminus\G'}(e) := \left\{
      \begin{array}{ll}
        e,       & \text{if $e\in s_\G\inv(V(\G'))$} \\
        s_\G(e), & \text{otherwise.}
      \end{array}\right.$
  }
\end{definition}
\begin{remark}\label{rmk:complement-pushout}
  There is a natural map $\G\setminus\G'\to\G$ which is injective on
  edges but not necessarily on vertices. It is easy to see that there
  is a pushout square~$\vcenter\inlinexy{{\G} &
    {\G\setminus\G'}\ar[l]\\{\G'}\ar[u] & \delta(\G\setminus\G')\ar[l]\ar[u]}$ in
  $\Graph$.  Topologically, $\G\setminus\G'$ is the complement in~$\G$
  of an open neighborhood of~$\G'$.
\end{remark}
Any other graph-related notions we use (such as \emph{tree},
\emph{geometric realization} and \emph{connected components}) will be
assumed well-known.

% We will  the usual graph notions of \emph{tree} and of
% \emph{connected components}.

%  \smallskip \begin{definition}
%    Given a graph $\G$, define its \emph{geometric realization} to be
%    the one-dimensional cell complex $|\G| := (V(\G)\disjcup
%    (H(\G)\cross[0,1]))/\!\!\sim$ where the equivalence relation is
%    generated by the identifications $(e,0)\sim s(e)$ and
%    $(e,t)\sim(r(e),1-t)$ for each $e\in H(\G)$, $t\in [0,1]$.
%  \end{definition}
%  Geometric realization defines a functor from $\Graph$ to topological
%  spaces and continuous maps. 

% \begin{definition}
%   Let $\G$ be a graph, and let $v, w\in V(\G)$. A \emph{path} joining
%   $v$ to $w$ is a sequence $e_1,\ldots,e_k$ where $e_i\in H(\G)$,
%   $v=s(e_1)$, $w=t(e_k)$, and $t(e_i)=s(e_{i+1})$ for
%   $i=0,\ldots,k-1$. The path is \emph{reduced} if $e_i\neq r(e_{i+1})$
%   for any $i=1,\ldots,k-1$. We say that $\G$ is a \emph{tree} if for
%   any pair $v, w\in V(\G)$, there is a unique reduced path joining $v$
%   to $w$.

%   Being connected by a path is an equivalence relation on vertices.
%   Each equivalence class of vertices, together with the set of
%   half-edges incident with them, forms a subgraph, which is called a
%   \emph{connected component}.  We will denote by $\pi_0(\G)$ the set
%   of connected components of $\G$.
% \end{definition}

\subsection{\Fgraphs}
\label{sec:fgraphs}
We will use the usual definition of \fgraph, except for the
restriction on vertices to be at least trivalent.
\begin{definition}\label{def:fgraph}
  A \emph{\fgraph} is a finite graph $\G$ equipped with a cyclic
  ordering on each of the sets $H(v)$, $v\in V(\G).$ We will encode
  this as a permutation $\sigma$ of $H(\G)$ with disjoint cycle
  decomposition given by the $H(v)$, and we will use the notation
  $\bsigma$ for the composition $\sigma\comp r: H(\G)\to H(\G)$, where
  $r$ is edge reversal. The \emph{boundary cycles} of $\G$ are the
  cycles of $\bsigma$.
\end{definition}

\smallskip\noindent It will be useful to represent the boundary of $\G$ as an
associated graph $\bdry\G$ having a natural morphism $\bdry\G\to\G$.
\begin{definition}
  Given a \fgraph $\G$, define $\bdry\G$ as follows.  Let the vertices
  of $\bdry\G$ be $V(\bdry\G) := H(\G)$. Let the half-edges be given
  by $H(\bdry\G) := H(\G)\cross\set{0,1}$, with edge-reversal
  involution $r(e,i) := (e,1-i)$. Define the attachment of half-edges
  by the source map $s(e,0) := e,$ $s(e,1) := \bsigma(e)$.  Define a
  morphism $\iota:\bdry\G\to\G$ by letting $\iota(e) := s(e)$ on
  vertices and $\iota(e,0) := e,$ $\iota(e,1) := r(e)$ on half-edges.
\end{definition}
It is immediate that the boundary cycles of $\G$ are in bijection with
the connected components of $\bdry\G$, each of which is a cyclic
graph. Moreover, $\bdry\G$ has a natural orientation induced by
$\bsigma$.

% \begin{proposition}
%   \label{prop:pushout}
%   Geometric realization takes pushout squares in $\Graph$ to pushout
%   squares of topological spaces.
% \end{proposition}

% --------------------------------------------------------------------
\subsection{\Ofgraphs}
% --------------------------------------------------------------------
Now we extend the notion of \fgraph to include free boundary with
labels in~$\B$.

\begin{definition}\label{def:ofgraph}
  An \ofgraph (with brane labels drawn from $\B$) is a \fgraph $\G$
  together with:
  \begin{enumerate}
  \item a distinguished \emph{free boundary} subgraph
    $\bdryF\G\subseteq\bdry\G$ such that the restriction of
    $\iota:\bdry\G\to\G$ is an embedding $\iotaF:\bdryF\G\inclto\G$,
    and,
  \item a labeling $\beta:\bdryF\G\to\B$ assigning an element of $\B$
    to each connected component of $\bdryF\G$. 
  \end{enumerate}
  Given an \ofgraph, we define its \emph{string boundary} $\bdrySt\G$
  as $\bdry\G\setminus\bdryF\G$
  (Definition~\ref{def:subgraph-complement}), and we let
  $\iotaSt:\bdrySt\G\to\G$ be the restriction of $\iota$.
\end{definition}
Each of $\bdryF\G$ and $\bdrySt\G$ is necessarily a disjoint union of
linear and cyclic graphs. In the case of $\bdryF\G$, the definition
includes the possibility of components that are isolated vertices.
Moreover, the two graphs $\bdryF\G$ and $\bdrySt\G$ intersect by
definition only on the set $\delta\bdrySt\G$ of vertices which are
endpoints of linear components of $\bdrySt\G$, and the restriction of
$\beta$ then makes $|\bdrySt\G|$ a $\B$-labeled one-manifold.

\begin{example}
  Consider the \fgraph\\[-1.5em]
  \begin{displaymath}
    \G = \picbox{mp/fbgsquare.4012},
  \end{displaymath}\\[-1em] 
  with the usual convention that the cyclic ordering at
  each vertex is counterclockwise.  Its boundary graph is\\[-1.5em]
  \begin{displaymath} \bdry\G = \picbox{mp/fbgsquare.4052}. \end{displaymath}\\[-1em]
  We can specify an \oc structure on $\G$ by choosing a subgraph
  $\bdryF\G$ of $\bdry\G$ with a $\B$-labeling of its components. This
  can be represented as follows:\\[-1.5em]
  \begin{displaymath} \picbox{mp/fbgsquare.4002}\!\!\!\!,
  \end{displaymath}\\[-1em] where the dotted lines and the hollow vertex are
  $\bdryF\G$ and the numbers indicate the labeling in
  $\B=\set{1,2,\ldots}$.  This gives $\G$ the structure of a genus one \ofgraph, having string boundary of type $I_{3,3}\disjcup
  I_{1,1}\disjcup S^1\disjcup S^1$.
\end{example}

\begin{example}
  Another example of an \ofgraph, using the same conventions, is\\[-1.5em]
  \begin{displaymath}
    \picbox{mp/fbgsquare.4001}.
  \end{displaymath}\\[-1em]
  With a boundary partitioning, this represents the ``composition''
  operation from the introduction.
\end{example}

\medskip As is well-known, given two graphs $\G_1, \G_2$ with $\G_1$ a
\fgraph, a morphism $\phi:\G_1\to\G_2$ induces a \fgraph structure on
$\G_2$ if it is \emph{simple}. A simple morphism is one that can be
written as a composition of ``edge collapses,'' that is, morphisms
whose effect on the $\B$-graph, up to isomorphism, is to collapse a
single non-loop edge down to a vertex, leaving the rest of the graph
intact.  More precisely:
\begin{definition}
  \label{def:simplemor} \label{def:ofg-mor} Call a morphism in
  $\Graph$ $\phi:\G_1\to\G_2$ \emph{simple} if the inverse image
  subgraph of any vertex of $\G_2$ is a tree and $\phi_H$ is injective
  on $\phi_H\inv(H(\G_2))$.  A morphism $\phi:\G_1\to\G_2$ of \fgraphs
  is a simple morphism of graphs which takes the \fgraph structure of
  $\G_1$ to that of $\G_2$.  

  \defsep Given \ofgraphs $\G_1$, $\G_2$, define a \emph{morphism}
  from $\G_1$ to $\G_2$ to be a morphism $\phi:\G_1\to\G_2$ of the
  underlying \fgraphs such that the induced morphism
  $\bdry\G_1\to\bdry\G_2$ restricts to a labeling-preserving simple
  morphism $\bdryF\G_1\to\bdryF\G_2$.  Let $\Fat_\B$ be the resulting
  category of \ofgraphs.
\end{definition}

% \begin{remark}
%   A morphism of \ofgraphs is a composition of edge collapses, but
%   there are restrictions on which edges may be collapsed.  A non-loop
%   edge $e$ may be collapsed only if the image of $\bdryF\G$ in $\G$
%   remains embedded afterwards.
% \end{remark}

\begin{definition}
  A \emph{partitioning} of an \ofgraph is a decomposition of
  $\bdrySt\G$ as a disjoint union of graphs $\bdry^-\G$
  (\emph{incoming})and $\bdry^+\G$ (\emph{outgoing}).  We denote by
  $\iota_\pm$ the restriction to $\bdry^\pm\G$ of the morphism
  $\iota:\bdry\G\to\G$.
\end{definition}

% --------------------------------------------------------------------
\section{$\B$-labeled graphs and their mapping spaces}
\label{sec:bgraphs}
% --------------------------------------------------------------------
Now we will let the set $\B$ index actual branes $\set{L_b}_{b\in\B}$
in a manifold $M$. Then, the $\B$-labels carried by an \ofgraph $\G$
(and by its boundary graph) may be used to specify constraints on maps
from $|\G|$ to $M$. This can be generalized as follows.

% Let $\G$ be an \ofgraph, so that it has an
% associated subgraph $B(\G)=\iota(\bdryF\G)\subseteq\G$ which carries a
% label in $\B$ on every connected component. This data can be used to
% specify constraints on maps from $\G$ to $M$.  

\begin{definition}
  A $\B$-labeled graph (or \bgraph for short) is a diagram
  \begin{displaymath}
    \G\xfrom{\theta} B(\G)\xto{\beta} \B
  \end{displaymath}
  of graphs in which $\theta$ is an embedding. The \bgraphs form a
  category $\Graph_\B$, in which a morphism $\phi:\G_1\to\G_2$ of
  $\B$-graphs is defined to be a pair $\phi:\G_1\to\G_2$,
  $\phi_B:B(\G_1)\to B(\G_2)$ of morphisms in $\Graph$ making the
  diagram
  \begin{displaymath}\vcenter\inlinexy{
      {B(\G_1)}\ar[rr]^{\beta_{1}}\ar[dd]_{\theta_{1}}\ar[dr]_{\phi_B} &                                                & {\B}  \\
      & {B(\G_2)}\ar[ur]_{\beta_{2}}\ar[dd]^{\theta_{2}} & {}    \\
      {\G_1}\ar[dr]_{\phi}                                        &                                                & {}    \\
      {}                                                        & {\G_2}                                          & {}    \\
    }\end{displaymath}  
  commute.
\end{definition}
This is another way of saying that a \bgraph is a graph $\G$ together
with a distinguished family $B(\G)$ of disjoint connected subgraphs,
each labeled by an element of $\B$. We put it in this way to clarify
the morphisms. We have two examples of \bgraphs.
\begin{example}
  Given an \ofgraph $\G$, and letting $B(\G)=\bdryF\G$, $\G$ is
  naturally a \bgraph via the diagram
  $\G\xinclfrom{\iotaF}\bdryF\G\xto{\beta}\B$.
\end{example}

\begin{example}
  The string boundary graph $\bdrySt\G$ of an \ofgraph becomes a
  \bgraph by letting $B(\bdrySt\G)$ be the discrete subgraph
  $\delta\bdrySt\G$ consisting of the endpoints of the linear
  components of $\bdrySt\G$.  There is a morphism
  $\iotaSt:\bdrySt\G\to\G$ of \bgraphs.
\end{example}

Given a \bgraph $\G$ and a family of closed submanifolds $L_b$ of a
given smooth closed manifold $M$, we may consider maps from $|\G|$ to
$M$ which respect the $\B$-labeling.
\begin{definition}[mapping spaces of \bgraphs]
  We will call $\M=(M,\set{L_b}_{b\in\B})$ a \emph{$\B$-brane system}
  in $M$.  Given a \bgraph $\G\xfrom{\theta} B(\G)\xto{\beta} \B$,
  define $\bmapm{\G}$ to be the space of continuous maps $f:|\G|\to M$
  such that for every connected component $C$ of $B(\G)$, the
  composition $f\comp|\theta|:|C|\to M$ has image in
  $L_{\beta(C)}\subseteq M$.
\end{definition}
For a fixed $\M$, this is a contravariant functor from $\Graph_\B$ to
topological spaces. We will denote its action on a morphism $\phi$ by
$\bmapm{\phi}$.  We will write $\bma{\blank}$ instead of
$\bmapm{\blank}$ when it is clear by context.

% --------------------------------------------------------------------
\section{Generalized Pontrjagin-Thom collapse}
\label{sec:thom-collapse}
% --------------------------------------------------------------------
The construction of the \oc string topology operations will
make use of a homology-level umkehr homomorphism for the restriction
map $\bmapm{\G}\xto{\bma{\iota_-}}\bmapm{\bdry^-\G}$, where $\G$ is a
partitioned \ofgraph. As is the case in the closed case, this
homomorphism will arise from a Pontrjagin-Thom collapse. Here, we
gather without proof some properties of a generalized form of the
Pontrjagin-Thom collapse. In the following section, we specialize to
the case of the maps that we are interested in.

\begin{theoremNoheader}[Basic fact]
  Assume given a homotopy pullback square 
$$\xymatrix{ {X}\ar[r]^f\ar[d]_p & {Y}\ar[d]^q \\
  {P}\ar[r]_g & {Q} }$$ such that $P$ and $Q$ are smooth closed
manifolds, and $g$ is a smooth map. There is a stable backwards map
$f^!:\Sigma^\infty Y_+\to X^{TQ-TP}$\!, where by abuse of notation
$X^{TQ-TP}$ stands for the Thom spectrum of the virtual bundle
$p^*g^*TQ-p^*TP$ on $X$.
\end{theoremNoheader}

If $h_*$ is a generalized homology theory and the manifolds $P$ and
$Q$ are \horiented, we may apply the Thom isomorphism theorem
$\tilde{h}_*(X^\xi)\xto{\isom} h_{*-\dim\xi}(X)$ to obtain a
homomorphism $h_*(Y)\to h_{*-d}(X)$, also denoted $f^!$, where $d=\dim
Q-\dim P$.

\subsection{Sketch of the construction}
This construction and its properties do not seem to be published in
this generality; however, an upcoming article by R.~Cohen and
J.~Klein~\cite{cohen:generalized-collapse} will include full
derivations.  We include only a sketch of the construction.

%The map $f^!$ is constructed as follows.  
First assume that the square is a pullback square with $g$ an
embedding and $q$ a locally trivial fiber bundle. In that case, we may
take a tubular neighborhood $U_g$ for $g$ and pull it back to an open
set $U_f\subseteq Y$. Then, $U_f$ will be a tubular neighborhood for
$f$, in the sense that the pair $(U_f,U_f\setminus f(X))$ is
homeomorphic to the pair $(p^*\nu_f,p^*\nu_f\setminus X)$, where
$\nu_f\downto P$ is the normal bundle of $Q$ and $X$ includes in
$p^*\nu_f$ as the zero section.  The desired map $f^!$ then comes
about as usual, by collapsing the complement of $U_f$ to a point and
mapping the rest homeomorphically.

If $g$ is not an embedding, then we choose an embedding $i$ of $P$
into a high-dimensional Euclidean space $\R^N$ and replace $g$ and $q$
respectively by $(g,i):P\to Q\cross \R^N$ and $q\cross\id: Y\cross
\R^N\to Q\cross\R^N$. Different choices of (sufficiently large) $N$
yield target Thom spaces that differ by a suspension, and the result
is a stable map into the Thom spectrum.

Finally, if $q$ is not locally trivial, it may be replaced by a Serre
fibration using the standard mapping path space construction (see,
e.g., \cite{may:concise}).  This fibration has sufficient structure to
allow a tubular neighborhood for $g$ to be lifted to one for $f$ using
parallel transport along paths.
% We refer to
% \cite{cohen:generalized-collapse} for the details.

\subsection{Naturality properties}
The generalized Pontrjagin-Thom collapse satisfies the following two
naturality properties.
\begin{proposition}[Functoriality]
  \label{prop:thomnat1} Consider a diagram
  \begin{displaymath}\xymatrix{ {X}\ar[r]^{f_1}\ar[d]_p & {Y}\ar[r]^{f_2}\ar[d]^q & {Z}\ar[d]^r\\
    {P}\ar[r]_{g_1} & {Q}\ar[r]_{g_2} & {R,} }\end{displaymath} where  the $g_i$ are
smooth maps between \horiented closed manifolds and
both squares are homotopy pullbacks.  

Then, the umkehr homomorphisms satisfy~$(f_2\comp f_1)^!  = f_1^!\comp
f_2^!.$
\end{proposition}

\begin{proposition}[Compatibility with induced maps on homology]
  \label{prop:thomnat2} Consider a commutative diagram
  
  \begin{displaymath}
  \xymatrix{ 
    {}
     & {P_1}\ar@{<-}[ld]_{p_1}\ar[rr]^{g_1}\ar'[d]^{s}[dd]
      & {}
       & {Q_1}\ar@{<-}[ld]^{q_1}\ar[dd] \\
    {X_1}\ar[dd]_{u}\ar[rr]_(0.33){f_1}
     & {}
      & {Y_1}\ar[dd]^(0.25){v}
       & {}    \\
    {}
     & {P_2}\ar@{<-}[ld]^{p_2}\ar'[r]^{g_2}[rr]
      & {}
       & {Q_2}\ar@{<-}[ld]^{q_2} \\
    {X_2}\ar[rr]_{f_2}
     & {}
      & {Y_2}
       & {}    \\
  }
  \end{displaymath}
  such that each $p_i$ is a homotopy pullback of $q_i$ via $g_i$,
  the $P_i$ and $Q_i$ are \horiented manifolds, and
  the virtual bundles $s^*(g_2^*TQ_2-TP_2)$ and $g_1^*TQ_1-TP_1$ are
  stably equivalent.  Then the diagram
  \begin{displaymath}
  \xymatrix{
    {h_*(X_1)}\ar[d]_{u_*} & {h_*(X_2)}\ar[d]^{v_*}\ar[l]_{f_1^!} \\
    {h_*(Y_1)} & {h_*(Y_2)}\ar[l]^{f_2^!}
  }
  \end{displaymath}
  commutes.
\end{proposition}

% --------------------------------------------------------------------
\section{Umkehr maps induced by morphisms of \bgraphs}
\label{sec:umkehr-maps-from-bgraph-morphisms}
% --------------------------------------------------------------------
We need a good family of morphisms $\phi:\G_1\to\G_2$ of \bgraphs for
which the induced map $\bma{\phi}:\bmapm{\G_2}\to\bmapm{\G_1}$ admits a
homology umkehr map.

First, observe that pushout squares in $\Graph_\B$ become pullback
squares of mapping spaces after applying $\bmapm{\blank}$. The proof
is standard and we omit it.
% \begin{proposition}\label{prop:respect-pushouts}
%   Let $\M$ be a $\B$-brane system, and let   
%   \begin{equation}
%     \label{diag:po}
%     \vcenter{\xymatrix{
%       {\G_1}&{\G_2}\ar[l]\\
%       {\Lambda_1}\ar[u]&{\Lambda_2}\ar[u]\ar[l] 
%     }}
%   \end{equation}
%   be a pushout square in $\Graph_\B$. Then, the induced diagram
%   \begin{equation}
%     \label{diag:pb}
%     \vcenter{\xymatrix{
%       {\bmapm{\G_1}}&{\bmapm{\G_2}}\ar@{<-}[l]\\
%       {\bmapm{\Lambda_1}}\ar@{<-}[u]&{\bmapm{\Lambda_2}}\ar@{<-}[u]\ar@{<-}[l] 
%     }}
%   \end{equation}
%   is a pullback square of topological spaces.
% \end{proposition}

% \smallskip 

Next, we identify a class of morphisms $\phi:\G_1\to\G_2$ in
$\Graph_\B$ with the property that the induced map $\bmapm{\phi}$
fibers naturally over a smooth map of manifolds.  For this, we
introduce the following construction.

\begin{definition} Given a \bgraph $\G$ define its \emph{vertex
$\B$-graph} $\V(\G)$ as the \bgraph that results from removing all the
edges of $\G$, keeping only the vertices and their labels.  Formally,
if $\G$ is given by a diagram $\G\xfrom{\theta} B(\G)\xto{\beta} \B$,
define $\V(\G)$ by the diagram 
  \begin{displaymath}V(\G)\xfrom{\restr{\theta}{V(B(\G))}} V(B(\G))
\xto{\restr{\beta}{V(B(\G))}} {\B}.\end{displaymath}
\end{definition} This defines a functor $\Graph_\B\to\Graph_\B$ having
a natural inclusion $\V(\G)\inclto\G$.

\begin{proposition}\label{prop:vert-pushout}
  Let $\phi:\G_1\to\G_2$ be a morphism in $\Graph_\B$ such that 
  \begin{enumerate}
  \item $\phi_H$ is a bijection $H(\G_1)\to H(\G_2)$, and,
  \item $\phi$ carries unlabeled edges to unlabeled edges; that is,
    $\phi$ takes edges not in the image of $B(\G_1)$ to edges not in
    the image of $B(\G_2)$.
  \end{enumerate}
  Then, the diagram
  $\vcenter\inlinexy{{\G_1}\ar[rr]^{\phi}&&{\G_2}\\{\V(\G_1)}\ar[u]\ar[rr]_{\V(\phi)}&&{\V(\G_2)}\ar[u]}$
  is a pushout square in $\Graph_\B$.
\end{proposition}
\begin{proof}
  The hypothesis says that, up to isomorphism compatible with $\phi$,
  $\G_2$ is obtained from $\G_1$ by identifying together vertices with
  a common preimage and then adjoining (possibly labeled) isolated
  vertices corresponding to $V(\G_2)\setminus \phi(V(\G_1))$. But this
  is equivalent to this square being a pushout.
\end{proof}

Note that $\bmapm{\V(\G)}\homeo \prod_{v\in V(\G)}L_v,$ where $L_v$ is
the submanifold of $M$ corresponding to the unique label carried by
the vertex $v$, or $M$ if $v$ is unlabeled. In particular,
$\bmapm{\V(\G)}$ is a smooth manifold.  In addition, if
$\phi:\G_1\to\G_2$ is a morphism in $\Graph_\B$, then $\phi$ induces a
smooth map $\bmapm{\V(\G_2)}\to\bmapm{\V(\G_1)}$; in fact, this map is
a cartesian product of coordinate projections of the form $L\cross
L'\projto L$ and diagonal inclusions of the forms $L\inclto L^p\cross
M^q$ and $M\inclto M^p$.

\medskip In view of this and Proposition~\ref{prop:vert-pushout}, we can make
the following definition.
\begin{definition}\label{def:easy-graph-umkehr}
  Let $\phi:\G_1\to\G_2$ be a morphism of \bgraphs satisfying the
  hypothesis of Proposition~\ref{prop:vert-pushout}.  Suppose given a
  $\B$-brane system $\M=\tuple{M,\set{L_b\subseteq M}_{b\in\B}}$ which
  is \horiented; that is, such that $M$ and each of the $L_b$ is
  oriented with respect to a multiplicative homology theory $h_*$.

  Define
  \begin{displaymath}\bmapm{\phi}^!:h_*(\bmapm{\G_1})\to h_*(\bmapm{\G_2})\end{displaymath} 
  as the umkehr homomorphism associated   by the generalized Pontrjagin-Thom collapse to the pullback square
  \begin{displaymath}\xymatrix{{\bmapm{\G_1}}\ar@{<-}[r]^{\bmapm{\phi}}&{\bmapm{\G_2}}\\{\bmapm{\V(\G_1)}}\ar@{<-}[u]\ar@{<-}[r]&{\bmapm{\V(\G_2)}.}\ar@{<-}[u]}\end{displaymath}
  This uses that the normal bundle of $L$ in $M$ is \horiented if both
  $TM$ and $TL$ are.
\end{definition}

\subsection{Enlarging the class of morphisms: the category  $\Graph_\B^!$}
\label{sec:more-umkehrables}

Thus a morphism of \bgraphs which is a bijection on edges and
preserves unlabeled edges induces an umkehr map in a natural way.
However, we will need umkehr homomorphisms for a larger class of
\bgraph morphisms:
\begin{definition}\label{def:Graph-bang}
  Let $\Graph_\B^!$ be the subcategory of $\Graph_\B$ consisting of
  morphisms $\phi:\G_1\to\G_2$ of \bgraphs such that:
  \begin{enumerate}
  \item $\phi_H$ induces an injection $H(\G_1)\to H(\G_2)$ of
    half-edges, and,
  \item $\phi$ carries unlabeled edges to unlabeled edges in the sense
    of Proposition~\ref{prop:vert-pushout}.
  \end{enumerate}
\end{definition}
We will construct the desired homomorphisms by showing that morphisms
in $\Graph_\B^!$ can be naturally factored up to homotopy into
morphisms satisfying the hypothesis of
Proposition~\ref{prop:vert-pushout}.

\medskip Let $\phi:\G_1\to\G_2$ be a morphism in $\Graph_\B^!$. Let
\begin{displaymath}\Xi_\phi=\bigdisjcup_{e\in  E(\G_2)\setminus\phi(E(\G_1))} \hat{e},\end{displaymath} 
where each $\hat{e}$ is a \bgraph consisting of two vertices joined by
a single edge, which is labeled by the same label carried by $e$ in
$\G_2$, or unlabeled if $e$ is unlabeled. Let $\xi_\phi$ be the
\bgraph obtained from $\Xi_\phi$ by collapsing each $\hat{e}$ to a
vertex carrying the same label, if any. Then, $\phi$ decomposes as
\begin{displaymath}
\G_1\xto{\phi'}\G_1\disjcup\xi_\phi\xfrom{\phi''}\G_1\disjcup\Xi_\phi\xto{\phi'''}\G_2,
\end{displaymath}
where $\phi'$ is the natural inclusion, $\phi''$ extends the defining
quotient map $\Xi_\phi\projto\xi_\phi$ by the identity on $\G_1$, and
$\phi'''$ is the morphism which extends $\phi$ by taking $\hat{e}$ to
$e$. 

The morphisms $\phi'$ and $\phi'''$ are readily seen to satisfy the
hypothesis of Proposition~\ref{prop:vert-pushout}. Moreover, $\phi''$
induces a homotopy equivalence
$\bmapm{\G\disjcup\xi_\phi}\to\bmapm{\G\disjcup\Xi_\phi}$, and
$\bmapm{\phi}$ is homotopic to
$\bmapm{\phi'}\comp\bmapm{\phi''}\inv\comp\bmapm{\phi'''}$.

\begin{example}
  We illustrate this for simplicity when $B(\G_1)=B(\G_2)=\emptyset$.
  Consider the morphism
  $\picbox{mp/fbgsquare.3014}\xto{\phi}\picbox{mp/fbgsquare.3015}$,
  which clearly lies in $\Graph_\B^!$. In this case, the factorization
  takes the form
  \begin{displaymath}
    \xymatrix{
    {\picbox{mp/fbgsquare.3014}} \inclar[r]^-{\phi'} &
    {\picbox{mp/fbgsquare.3017}} & 
    {\picbox{mp/fbgsquare.3016}} \surjar[l]_-{\phi''} \surjar[r]^-{\phi'''} &
    {\picbox{mp/fbgsquare.3015}.}}
\end{displaymath}
\end{example}

\medskip\noindent We may now make the following definition.
\begin{definition}\label{def:graph-umkehr}
  Given a morphism $\phi\in\Graph_\B^!$ and an \horiented
  $\B$-brane system $\M$ we define the \emph{umkehr homomorphism}
  \begin{displaymath}
  \bmapm{\phi}^!:h_*(\bmapm{\G_1})\to h_*(\bmapm{\G_2})
  \end{displaymath}
  as the composition $\bmapm{\phi'''}^!\comp (\bmapm{\phi''})_*\comp
  \bmapm{\phi'}^!$, where the homomorphisms $\bmapm{\phi'}^!$ and
  $\bmapm{\phi'''}^!$ are given by
  Definition~\ref{def:easy-graph-umkehr}.
\end{definition}
\begin{remark}
  Since $\xi_\phi$ is discrete, the map $\bmapm{\phi'}$ is a
  projection $\bmapm{\G_1}\cross N \to \bmapm{\G_1}$ with $N$ a closed
  \horiented manifold. Its corresponding umkehr map is simply crossing
  with the fundamental class of $N$. In these terms, the map
  $\bmapm{\phi''}$ is the inclusion $\bmapm{\G_1}\cross N\to
  \bmapm{\G_1}\cross PN$, where $PN$ stands for the space of arbitrary
  continuous paths in $N$, with $N$ included as the constant paths.
\end{remark}

Finally, we observe that this homomorphism behaves well under the
appropriate notion of simple morphism for \bgraphs, as well as under
pushouts of \bgraph embeddings.
\begin{definition}
  Say that a morphism $\phi=(\phi,\phi_B):\G_1\to\G_2$ of \bgraphs is
  \emph{simple} if each of $\phi:\G_1\to\G_2$ and $\phi_B:B(\G_1)\to
  B(\G_2)$ is a simple morphism of graphs.
\end{definition}
A morphism of \ofgraphs is in particular a simple morphism of the
underlying \bgraphs, and it induces a simple morphism of the
associated string boundary \bgraphs.

\begin{remark}\label{rmk:wholly-monochromatic}
  A simple morphism of \bgraphs is of course a composition of edge
  collapses.  However, a non-loop edge may be collapsed only if it is
  ``wholly monochromatic,'' that is, if the subgraph consisting of the
  edge and its two endpoints is either disjoint from $B(\G)$ or
  contained in $B(\G)$.
\end{remark}

\begin{proposition}
  \label{prop:umkehr-simple} 
  Let
  $\vcenter{\inlinexy{{\G_1}\ar[r]^{\alpha}\ar[d]_{\gamma_1}&{\G_2}\ar[d]^{\gamma_2}\\{\tG_1}\ar[r]_{\beta}&{\tG_2}}}$
  be a commutative diagram in $\Graph_\B$. Suppose that $\alpha$ and
  $\beta$ lie in $\Graph_\B^!$, and that either:
  \begin{enumerate}
  \item the diagram is a pushout and $\alpha$ is an embedding, or,
  \item the $\gamma_i$ are simple.
  \end{enumerate}
  Then, the diagram
  \begin{equation} \label{diag:hsquare}\xymatrix{{h_*(\bmapm{\G_1})}\ar[r]^{\bma{\alpha}^!}&{h_*(\bmapm{\G_2})}\\{h_*(\bmapm{\tG_1})}\ar[r]_{\bma{\beta}^!}\ar[u]^{\bma{\gamma_1}_*}&{h_*(\bmapm{\tG_2})}\ar[u]_{\bma{\gamma_2}_*}}
  \end{equation}
  commutes.
\end{proposition}
\begin{proof}
  \emph{Case 1:} the diagram is a pushout and $\alpha$ is an
  embedding. In this case, the $\gamma_i$ induce a commutative diagram
  \begin{displaymath}\xymatrix{
      {\G_1}\ar[d]\ar[r]^{\alpha'} & {\G_1\disjcup\xi}\ar[d] & {\G_1\disjcup\Xi}\ar[d]\ar[l]_{\alpha''}\ar[r]^{\alpha'''} & {\G_2}\ar[d]  \\
      {\tG_1}\ar[r]^{\beta'} & {\tG_1\disjcup\xi} & {\tG_1\disjcup\Xi}\ar[l]_{\beta''}\ar[r]^{\beta'''} & {\tG_2,}  \\
    }
  \end{displaymath}
  where $\xi=\xi_\alpha\isom\xi_\beta$,
  $\Xi=\Xi_\alpha\isom\Xi_\beta$, and in which the leftmost and
  rightmost squares are pushouts in $\Graph_\B$. It follows that these
  squares become pullbacks upon applying $\bmapm{\blank}$. Moreover,
  each of these pullbacks fibers, in the sense of
  Proposition~\ref{prop:thomnat1}, over the corresponding pullback
  square of manifolds obtained by applying $\bmapm{\V(\blank)}$. The
  result follows because the hypothesis on normal bundles of
  Proposition~\ref{prop:thomnat1} is easily verified to hold for the
  latter squares.

  \smallskip \newcommand{\he}{\hat{e}} \emph{Case 2:}  $\gamma_1$
  and $\gamma_2$ are simple. By induction we can assume that
  $\gamma_2$ collapses a single edge.  If the collapsed edge is in the
  image of $\alpha$, the result reduces to case~1.  Assume then that
  $\gamma_2$ collapses a single edge $e$ which is not in the image of
  $\alpha$.  We have up to isomorphism that $\tG_2 = \G_2/e$ and
  $\tG_1 = \G_1$.  Let us change notation for clarity, writing $\Xi$
  and $\xi$ for $\Xi_\beta,$ $\xi_\beta$ respectively. We clearly have
  $\Xi_{\alpha} \isom \Xi\disjcup \hat{e}$ and
  $\xi_\alpha\isom\xi\disjcup\pt_e$, where $\pt_e$ stands for a
  one-vertex graph carrying the same label as $e$, if any. Also write
  $\gamma$ for $\gamma_2$; $\gamma_1$ becomes the identity in this
  case.
  
  We have the diagram
  \begin{displaymath}
    \xymatrix{
      \G_1\ar[r]^-{\alpha'}\ar[ddr]_{\beta'} & \G_1\disjcup\xi\disjcup\pt_e & \G_1\disjcup\Xi\disjcup\he\ar[l]_{\alpha''}\ar[r]^-{\alpha'''} & \G_2\ar[dd]^{\gamma} \\
      \\
         &  \G_1\disjcup\xi           & \G_1\disjcup\Xi\ar[l]^{\beta''}\ar[r]_{\beta'''}            & \G_2/e
    }
  \end{displaymath}
  and we are to show that
  $\bma{\alpha'''}^!\comp\bma{\alpha''}_*\comp\bma{\alpha'}^!=\bma{\gamma}_*\comp\bma{\beta'''}^!\comp\bma{\beta''}_*\comp\bma{\beta'}^!.$
  For this, complete the diagram as follows:
  \begin{displaymath}
    \xymatrix{
      \G_1\ar[r]^-{\alpha'}\ar[ddr]_{\beta'} & \G_1\disjcup\xi\disjcup\pt_e & \G_1\disjcup\Xi\disjcup\he\ar[l]_{\alpha''}\ar[r]^-{\alpha'''}\ar[d]^{\gamma'} & \G_2\ar[dd]^{\gamma} \\
                                          &                          & \G_1\disjcup\Xi\disjcup\pt_e\ar[ul]^{\phi}\ar[dr]^{\psi}                 \\
                                          &  \G_1\disjcup\xi\ar[uu]^{\delta}    & \G_1\disjcup\Xi\ar[l]^{\beta''}\ar[r]_{\beta'''}\ar[u]^{\delta'}              & \G_2/e
    }
  \end{displaymath}
  Here, $\gamma'$~collapses~$\he$ to~$\pt_e$, $\phi$ extends~$\beta''$
  by the identity on~$\pt_e$, $\psi$~extends~$\beta'''$ by mapping~$\pt_e$
  to the vertex to which $e$ is collapsed, and $\delta$ and $\delta'$
  are the obvious inclusions.  Note that $\delta$, $\delta'$ and
  $\psi$ satisfy the hypothesis of
  Definition~\ref{def:easy-graph-umkehr}.

  By Remark~\ref{rmk:wholly-monochromatic} above, $e$ is wholly
  monochromatic; this ensures that $\bma{\delta}$ and $\bma{\delta'}$
  have the same normal bundle data (in the sense of
  Proposition~\ref{prop:thomnat2}) when $\bma{\blank}$ is applied (the
  stable normal bundle is a pullback of $-TL$ for both, where $L$ is
  either $M$ or the brane submanifold corresponding to the label
  carried by $e$).  This is also true for the pair $\bma{\alpha'''}$,
  $\bma{\psi}$.

  The equality of the two homomorphisms then follows by applying the
  naturality properties:
  \begin{eqnarray*}
    \bma{\gamma}_*\comp\bma{\beta'''}^!\comp\bma{\beta''}_*\comp\bma{\beta'}^! 
    & = & (\bma{\gamma}_*\comp \bma{\psi}^!) \comp \bma{\delta'}^! \comp\bma{\beta''}_*\comp\bma{\beta'}^! \\
    & = & \bma{\alpha'''}^!\comp \bma{\gamma'}_* \comp (\bma{\delta'}^! \comp\bma{\beta''}_*)\comp\bma{\beta'}^! \\
    & = & \bma{\alpha'''}^!\comp \bma{\gamma'}_* \comp \bma{\phi}_* \comp(\bma{\delta}^!\comp\bma{\beta'}^!) \\
    & = & \bma{\alpha'''}^!\comp (\bma{\gamma'}_* \comp \bma{\phi}_*) \comp\bma{\alpha'}^! \\
    & = & \bma{\alpha'''}^!\comp \bma{\alpha''}_* \comp\bma{\alpha'}^!, \\
  \end{eqnarray*}
  as desired. 
\end{proof}

% --------------------------------------------------------------------
\section{Definition of the operations}
\label{sec:operations}
% --------------------------------------------------------------------

With the constructions of the previous section, the definition of the
string topology operation associated to an \ofgraph is
straightforward.
\begin{definition}\label{def:admissible-partitioning}
  Say that a partitioning of an \ofgraph $\G$ is \emph{admissible} if
  the morphism $\iota_-:\bdry^-\G\to\G$ lies in $\Graph_\B^!$. If this
  is the case, we say that $\G$ is \emph{well-partitioned}.

  Define a category $\FatP_\B$ as follows. An object of $\FatP_\B$ is
  a well-partitioned \ofgraph.  The morphisms from $\G_1$ to $\G_2$
  are the morphisms $\phi:\G_1\to\G_2$ of \ofgraphs which respect the
  partitioning, in the sense that the induced morphism
  $\bdrySt\G_1\to\bdrySt\G_2$ takes~$\bdry^-\G_1$ to~$\bdry^-\G_2$ and
  $\bdry^+\G_1$~to~$\bdry^+\G_2$.
\end{definition}

\begin{definition}\label{def:operations}
  Let $\M$ be an \horiented $\B$-brane system, and let $\G$ be a
  well-\linebreak partitioned \ofgraph.  Define the homomorphism $\G_*$ as
  \begin{displaymath}
    \G_*:= \bma{\iota_+}_*\comp\bma{\iota_-}^!: h_*(\bmapm{\bdry^-\G}) \to h_{*}(\bmapm{\bdry^+\G}). \\
  \end{displaymath}
  This homomorphism is the \emph{\oc string topology
    operation} corresponding to~$\G$.
\end{definition}
The operations $\Gamma_*$ are invariant under morphisms of \ofgraphs:
\begin{proposition}
  \label{prop:collapse-ind}
  If $\G_1,\G_2\in\FatP_\B$ and there is a morphism
  $\phi:\G_1\to\G_2$, then $(\G_2)_* =
  \bma{\bdry^+\phi}_*\inv\comp(\G_1)_*\comp\bma{\bdry^-\phi}_*$.
\end{proposition}
\begin{proof}
  In the commutative diagram
  \begin{displaymath}
  \xymatrix{
    {\bdry^+\G_1}\ar[r]^{\bdry^+\phi}\ar[d]_{\iota^1_+} & {\bdry^+\G_2}\ar[d]^{\iota^2_+}\\
    {\G_1}\ar[r]^{\phi}                                 & {\G_2}\\
    {\bdry^-\G_1}\ar[u]^{\iota^1_-}\ar[r]^{\bdry^-\phi} & {\bdry^-\G_2,}\ar[u]_{\iota^2_-}
  }
  \end{displaymath}
  of \bgraphs, the morphisms $\iota^{1}_-, \iota^{2}_-$ lie in
  $\Graph^!_\B$ and the horizontal morphisms are simple. Then, with an
  application of Proposition~\ref{prop:umkehr-simple}, we have that
  \begin{displaymath}\begin{array}{rcl}
    (\G_2)_*
    & = & \bma{\iota^2_+}_*\comp\bma{\iota^2_-}^! \\
    & = & (\bma{\bdry^+\phi}_*\inv\comp\bma{\iota^1_+}_*\comp\bma{\phi}_*)\comp(\bma{\phi}_*\inv\comp\bma{\iota^1_-}^!\comp\bma{\bdry^-\phi}_*) \\
    & = & \bma{\bdry^+\phi}_*\inv\comp\bma{\iota^1_+}_*\comp\bma{\iota^1_-}^!\comp\bma{\bdry^-\phi}_* \\
    & = & \bma{\bdry^+\phi}_*\inv\comp(\G_1)_*\comp\bma{\bdry^-\phi}_*, \\
  \end{array}\end{displaymath}
as desired
\end{proof}

\section{Gluing}
\label{sec:gluing}
% --------------------------------------------------------------------

Now we describe the combinatorial counterpart to gluing of cobordisms.

Suppose given $\G_1, \G_2\in\FatP_\B$ together with isomorphisms
$\bdry^+\G_1\xfrom{\gamma_1}\Delta\xto{\gamma_2}\bdry^-\G_2$, where
$\Delta$ is a \bgraph. Assume that $\gamma_2\comp\gamma_1\inv$ is
orientation-reversing. We may construct a \bgraph $\Gg$ by identifying
the outgoing boundary of $\G_1$ with the incoming boundary of $\G_2$
according to their common identification with $\Delta$.  More
precisely, we can define $\Gg$ by the pushout diagram
\begin{equation}\label{diag:gluepushout}\xymatrix{
    {\Gg}   & {\G_2}\ar[l] \\
    {\G_1}\ar[u] & {\Delta,}\ar[l]^{\alpha_1}\ar[u]_{\alpha_2} \\
  }\end{equation}
in~$\Graph_\B$, where $\alpha_1=\iota_-^1\comp\gamma_1$ and
$\alpha_2=\iota_+^2\comp\gamma_2$.

While the pushout $\Gg$ exists, it does not necessarily inherit an
\ofgraph structure having the correct isomorphism type. For that, we
need an extra condition on the partitioning:
\begin{definition}\label{def:very-admissible}
  Say that a partitioning of an \ofgraph $\G\in\FatP_\B$ is \emph{very
    admissible} if the inclusion $\iota_-:\bdry^-\G\to\G$ is an
  embedding of \bgraphs (in that case, $\G$ is \emph{very
    well-partitioned}).
\end{definition}
\begin{remark}
  This condition is somewhat analogous to the chord diagram constraint
  of \cite{cohen-godin:polarized-view-string-topo}. However, we don't
  require the complement of the incoming boundary to be a forest; we
  may do this because the factorization described in
  Section~\ref{sec:more-umkehrables} makes it unnecessary to collapse
  this complement.
\end{remark}

\begin{lemma}
  Suppose given $\G_1,\G_2\in\FatP_\B$, with $\G_2$ very
  well-partitioned. The \bgraph $\Gg$, as defined by pushout
  diagram~\eqref{diag:gluepushout}, has an induced partitioned
  \ofgraph structure with $\bdry^-\G_1\isom\bdry^-(\Gg)$,
  $\bdry^+\G_2\isom\bdry^+(\Gg)$ as \bgraphs.
\end{lemma}
\begin{proof}
  To describe the \fgraph structure permutation $\sigma$ of $\Gg$,
  we may first describe its associated permutation $\bsigma$
  (Definition~\ref{def:fgraph}), and then verify that, if we define
  $\sigma$ as $\bsigma\comp r$, the result is a \fgraph.

  Describing $\bsigma$ is equivalent to constructing a graph $\bdry\G$
  and a morphism $\iota:\bdry\G\to\G$ in $\Graph$, taking edges to
  edges, such that:
   \begin{enumerate}
   \item $\bdry\G$ is a disjoint union of cyclic graphs, each with a
     chosen orientation, and,
   \item every oriented edge $e\in H(\G)$ is the orientation-preserving
     image of exactly one edge in $\bdry\G$.
   \end{enumerate}
   In the process, we will construct a $\B$-labeled free boundary
   subgraph $\bdryF(\Gg)\subseteq\bdry(\Gg)$ and the partitioning
   $\bdrySt(\Gg)=\bdry^-(\Gg)\disjcup\bdry^+(\Gg)$.

   Since pushouts in $\Graph_\B$ can be constructed setwise, we have
   that $B(\Gg)= B(\G_1)\cup_{B(\Delta)} B(\G_2)$. The graph $B(\Gg)$
   is a disjoint union of cyclic and linear (possibly degenerate)
   graphs because each of the $B(\G_i)$ is by definition one such
   graph, and the discrete graph $B(\Delta)$ identifies pairs of
   endpoints of linear components. Thus we may identify the discrete
   graph $\delta B(\Gg)$ as consisting of the endpoints of the linear
   components of $B(\Gg)$, with multiplicity two for the degenerate
   ones. Define $\bdry(\Gg)$ by the pushout
   \begin{displaymath}\xymatrix{
       {\bdry(\Gg)} & {\bdry^-\G_1\disjcup\bdry^+\G_2}\ar[l] \\
       {B(\Gg)}\ar[u] & {\delta B(\Gg).}\ar[l]\ar[u] \\
     }
   \end{displaymath}
   (this is a case of the pushout in
   Remark~\ref{rmk:complement-pushout}).  It is easy to see that the
   resulting graph is a disjoint union of cyclic graphs. It inherits a
   preferred orientation from the orientations on the $\bdry\G_i$,
   contains $B(\Gg)$ as a subgraph, and satisfies $\bdry(\Gg)\setminus
   B(\Gg)\isom \bdry^-\G_1\disjcup\bdry^+\G_2$.

   It remains to verify that the permutation $\sigma:=\bsigma\comp r$
   implicit in this construction has exactly one disjoint cycle $H(v)$
   per vertex of $\Gg$. Here, we use the fact that $\G_2$ is very well
   partitioned.  This implies that the map $V(\G_1)\to V(\Gg)$ is
   injective, and hence that the vertices of the pushout $\Gg$ are of
   only two types:
   \begin{enumerate}
   \item vertices $u$ corresponding to those $v\in V(\G_2)$ not in the
     image of $V(\Delta)$, and,
   \item vertices $u$ corresponding to $w\in V(\G_1)$, resulting as
     the identification of $w$ with each of the vertices
     $\alpha_2(\alpha_1\inv(w))\subseteq V(\G_2)$.
   \end{enumerate}
   In the first case, we have $H(u)\isom H(v)$, and the cyclic
   ordering is given by the cyclic ordering in $H(v)$. 

   In the second case, for each $v\in\alpha_2(\alpha_1\inv(w))$,
   $H(v)$ has a preferred linear ordering given by opening its cyclic
   ordering at the unique $e\in H(v)$ lying in $\bdry^-(\G_2)$; it is
   unique since otherwise $\bdry^-\G_2\to\G_2$ would not be injective.
   Then, $H(u)$ is obtained from $H(w)$ by inserting the half-edges
   $H(v)$ in this linear order in spaces between half-edges in $H(w)$
   determined by the image of $\Delta$ in $\G_1$. This yields a cyclic
   ordering on $H(u)$ given by ``cyclically splicing'' the linear
   orders on the $H(v)$ into the cyclic ordering on $H(w)$.

   These cyclic orders are directly seen to agree with the ones
   induced by $\bdry(\Gg)$ above.
%
%    (FIXME---THIS IS VERY EASY TO SEE BUT HARD TO WRITE IN DETAIL; HOW
%    MUCH SPACE SHOULD I DEVOTE TO IT? SOME PICTURES OF WHAT HAPPENS AT
%    VERTICES WHEN GLUING MIGHT HELP.)
\end{proof}
\begin{remark}
  Verifying that $\sigma$ gives a cyclic ordering to each $H(u)$,
  $u\in V(\Gg)$ is a necessary step. If $\G_2$ is not \emph{very} admissible,
  then, while we may still construct the requisite $\bdry(\Gg)\to\Gg$,
  it may induce a permutation $\sigma$ for which a single $H(u)$
  contains multiple cycles of $\sigma$.
\end{remark}

\subsection{Compatibility with gluing}
\label{sec:gluing-compat}
We will show that the \oc string operations are compatible
with the gluing construction, in the sense that
$(\Gg)_*=(\G_2)_*\comp(\G_1)_*$, appropriately interpreted,
whenever the $\G_i$ are very well-partitioned \ofgraphs that fit in a
gluing setting.

\begin{proposition}
  \label{prop:gluing}
  Let gluing data
  $\tuple{\G_1,\G_2,\bdry^+\G_1\xfrom{\gamma_1}\Delta\xto{\gamma_2}\bdry^-\G_2}$
  be given with $\G_2$ very well-partitioned, and let be formed
  accordingly. Then,
  \begin{displaymath}
  (\Gg)_*=(\G_2)_*\comp\bma{\gamma_1\inv\comp\gamma_2}_*\comp(\G_1)_*.
  \end{displaymath}
\end{proposition}
\begin{proof}
  We have the diagram
  \begin{displaymath}\xymatrix{
    {\bdry^-\G_1}\ar[r]^{\iota_-^1}\ar[dr]_{\iota_-} & {\G_1}\inclar[d]^{\delta_1} & \Delta\ar[l]_{\alpha_1}\inclar[d]^{\alpha_2} \\
                                                     & {\Gg}                       & {\G_2}\ar[l]_{\delta_2} \\
                                                     &                             & {\bdry^+\G_2.}\ar[ul]^{\iota_+}\ar[u]_{\iota_+^2}
%   \ar@{} "1,2";"2,3" |{(\star)} 
  }\end{displaymath} 
where the upper-right square is the defining pushout.
Proposition~\ref{prop:umkehr-simple} (case 1) implies that this square
has the property
\begin{equation}
  \bma{\alpha_2}^!\comp\bma{\alpha_1}_* = \bma{\delta_2}_*\comp\bma{\delta_1}^!\label{eq:6}
\end{equation}
in homology.  Then we have:
\begin{displaymath}\begin{array}{rclr}
  (\Gg)_*
  & = & \bma{\iota_+}_*\comp\bma{\iota_-}^!\\%&\text{(definition)} \\
  & = & \bma{\iota^2_+}_*\comp\bma{\delta_2}_*\comp\bma{\delta_1}^!\comp\bma{\iota^1_-}^!\\%&\text{(functoriality of the umkehr construction and of $h_*$)} \\
  & = & \bma{\iota^2_+}_*\comp\bma{\alpha_2}^!\comp\bma{\alpha_1}_*\comp\bma{\iota^1_-}^!\\%&\text{(equation~\eqref{eq:6})}\\
  & = & \bma{\iota^2_+}_*\comp\bma{\iota^2_-\comp\gamma_2}^!\comp\bma{\iota^1_+\comp\gamma_1}_*\comp\bma{\iota^1_-}^!\\
  & = & \bma{\iota^2_+}_*\comp\bma{\iota^2_-}^!\comp\bma{\gamma_2}^!\comp\bma{\gamma_1}_*\comp\bma{\iota^1_+}_*\comp\bma{\iota^1_-}^!\\
  & = & (\G_2)_*\comp\bma{\gamma_1\inv\comp\gamma_2}_*\comp(\G_1)_*.
\end{array}\end{displaymath}
\end{proof}

% --------------------------------------------------------------------
\section{The string topology \btqft}
% --------------------------------------------------------------------
Let $h_*$ be a multiplicative homology theory whose coefficient ring
$h_*$ is a \emph{graded field}, that is, a graded ring in which all
nonzero homogeneous elements are invertible.  Given an
\horiented $\B$-brane system
$\M=(M,\set{L_b}_{b\in\B})$, we have a \bfam
\begin{displaymath}\fV_\M=(h_*LM,\set{h_*\P{M}{L_a,L_b}}_{a,b\in\B})\end{displaymath} over the
coefficient ring of $h_*$. Here $LM$ is the free loop space of $M$.
The constraints on $h_*$ have the effect that there is a product map
$h_*(X)\tensor_{h_*} h_*(Y)\to h_*(X\cross Y)$ (because $h_*$ is
multiplicative) which is moreover an isomorphism of graded
$h_*$-modules (because the graded field condition on $h_*$ makes the
K\"unneth spectral sequence collapse).

In this section, we will describe the positive-boundary \btqft
structure on $\fV_\M$ arising from the \oc string topology
operations from Definition~\ref{def:operations}.
\begin{definition}
  Recall that, given a \fgraph $\G$, there is an associated oriented
  surface $S(\G)$ having $|\G|$ as a deformation retract, and having
  an identification $\bdry S\isom|\bdry\G|$. If $\G$ is additionally
  an \ofgraph, $S$ becomes naturally an \osurf by decreeing the image
  of $|\bdryF\G|$ in $\bdry S$ to be the free
  boundary,\footnote{Strictly speaking we take a small closed
    neighborhood of $|\bdryF\G|$ in $\bdry S$ in order to deal with
    degenerate components of $\bdryF\G$.} with the induced
  $\B$-labeling. A partitioning of $\G$ then makes $S(\G)$ into an
  \ocobo.
\end{definition}

\begin{proposition}\label{prop:fattening-good}
  \begin{itemize}
  \item If two partitioned \ofgraphs are related by a morphism in
    $\FatP_\B,$ then their corresponding \ocobos are isomorphic

  \item Gluing of partitioned \ofgraphs translates, up to isomorphism,
    into gluing of the corresponding \ocobos.
  \end{itemize}
\end{proposition}
\begin{proofsketch}
  The first statement is clear from the corresponding result for
  \fgraphs. For the second statement, we may observe that the defining
  pushout diagram \eqref{diag:gluepushout} for gluing realizes to a
  homotopy pushout square equivalent to the counterpart diagram that
  arises when gluing cobordisms; this determines the Euler
  characteristic of $S(\Gg)$. By a comparison of boundary
  components, also the genus and $\B$-labeling structure are seen to
  correspond.
\end{proofsketch}
We will make use of the following lemma, whose proof we defer to the appendix.
\begin{lemma}\label{lem:components}
  The connected components of $\FatP_\B$ are in one-to-one
  correspondence with isomorphism types of \oc cobordisms~$S$
  for which $\bdry^+S$ intersects every connected component of~$S$.
\end{lemma}

Using the notation of Definition~\ref{def:pcobB}, we first show how a
morphism $S:x_-\to x_+$ in $\pCobpb_\B$ (with $x_\pm\in\mM_\B$)
produces a homomorphism $\sop_S:\fV_\M(x_-)\to \fV_\M(x_+)$.  We may
assume that $S$ is connected, and induce the remaining operations by
tensor product.

By Lemma~\ref{lem:components}, we can represent the morphism $S$ by a
well-partitioned \ofgraph $\G\in\FatP_\B$ having $|\bdry^-\G|\isom
|x_-|$ and $|\bdry^+\G|\isom |x_+|^*$ as $\B$-labeled one-manifolds,
together with:
\begin{enumerate}
\item a choice of orientation-preserving parametrization of each
  component of $\bdry\G$, and,
\item a choice $L_-$ of linear orderings of the incoming boundary
  components in each $\B$-labeling type, and a similar choice of $L_+$
  for the outgoing boundary components.
\end{enumerate}
A boundary parametrization may be specified by a starting point on
each closed string boundary component; it is then determined by
starting at that point and parameterizing each edge in a
piecewise-linear fashion at constant speed. String boundary intervals
have a natural parametrization in the same way.

The parametrizations and linear orderings (together with our
restrictions on $h_*$) give isomorphisms $h_*(\bdry^\pm\G)\isom
\fV_\M(x_\pm)$, and therefore the operation $\G_*$ defines a
homomorphism $\fV_\M(x_-)\to \fV_\M(x_+)$, which we call $\sop_S$.

\begin{lemma}
  The homomorphism $\sop_S$ is well-defined; that is, it is
  independent of the choice of representing well-partitioned \ofgraph
  and the choice of boundary parameterization.
\end{lemma}
\begin{proof}
  We will temporarily use $\sop_\G$ for the operation defined with a
  particular choice of boundary-parametrized \ofgraph $\G$.  Define a
  category $m\FatP_\B$ in which an object is a well-partitioned
  \ofgraph together with a choice of starting point in each string
  boundary cycle, with morphisms being those morphisms of partitioned
  \ofgraphs which preserve the choices of starting points. We can
  argue as in \cite{cohen-godin:polarized-view-string-topo} to show
  that the obvious forgetful functor $m\FatP_\B\to\FatP_\B$ is a torus
  fibration on each component, and thus the connected components of
  $m\FatP_\B$ are also in one-to-one correspondence with isomorphism
  types of \ocobos. Thus, it is enough to show that if
  $\phi:\G_1\to\G_2$ is a morphism in $m\FatP_\B$ with $\G_1$ and
  $\G_2$, then $\sop_{\G_1}=\sop_{\G_2}$.

  Consider the diagrams
  \begin{displaymath}\xymatrix{
    {} &
     {|\bdry^-\G_1|}\ar[dd]^{|\bdry^-\phi|} &
     {} &
        {|\bdry^+\G_1|}\ar[dd]_{|\bdry^+\phi|} & 
         {} \\
    {|x_-|}\ar[ur]^{\gamma_1^-}_{\isom}\ar[dr]_{\gamma_2^-}^{\isom} &
     &
      &
       &
        {|x_+|}\ar[ul]_{\gamma_1^+}^{\isom}\ar[dl]^{\gamma_2^+}_{\isom}\\
    {} &
     {|\bdry^-\G_2|} &
     {} &
        {|\bdry^+\G_2|.} & 
         {} \\
  }\end{displaymath}
  The $\gamma_i^\pm$ are isomorphisms of $\B$-labeled one-manifolds
  induced by the choices of starting points on the boundary cycles.
  The two triangles commute up to homotopy relative to the boundary.
  It follows that, when we apply the functor $\bmapm{\blank}$, the
  induced homomorphisms on homology satisfy
  \begin{displaymath}\bma{\gamma_2^\pm}_*=\bma{\gamma_1^\pm}_*\comp\bma{\bdry^\pm\phi}_*.\end{displaymath}
  By definition, we have
  \begin{displaymath}\sop_{\G_2} = \bma{\gamma_2^+}_*\comp(\G_2)_*\comp \bma{\gamma_2^-}_*\inv,\end{displaymath}
  and therefore 
  \begin{displaymath}\sop_{\G_2} = \bma{\gamma_1^+}_*\comp\bma{\bdry^+\phi}_*\comp(\G_2)_*\comp \bma{\bdry^-\phi}_*\inv\comp\bma{\gamma_1^-}_*\inv.\end{displaymath}
  By Proposition~\ref{prop:collapse-ind}, we have that $(\G_1)_* =
  \bma{\bdry^+\phi}_*\comp(\G_2)_*\comp \bma{\bdry^-\phi}_*\inv$, and thus
  \begin{displaymath}\sop_{\G_2} = \bma{\gamma_1^+}_*\comp(\G_1)_*\comp\bma{\gamma_1^-}_*\inv.\end{displaymath}
  But this is equal to $\sop_{\G_1}$ by definition.
\end{proof}

\begin{lemma}
  The assignment $S\mapsto \sop_S$ is functorial.
\end{lemma}
\begin{proof}
  This is directly implied by Proposition~\ref{prop:gluing}.
\end{proof}
\bigskip\noindent%
As a corollary, we have our theorem.
\begin{maintheorem}\label{thm:main1}
  If $h_*$ is a multiplicative generalized homology theory for which
  the coefficient ring is a graded field, then, letting $\fV_\M$ be
  the \bfam
  \begin{displaymath}\tuple{h_*LM,\set{h_*\P{M}{L_a,L_b}}_{a,b\in\B}}\end{displaymath}
  over $h_*(\pt)$, there is a positive-boundary \btqft structure on
  $\fV_\M$ which extends the known positive-boundary string topology
  \tqft structure on $h_*(LM)$.
\end{maintheorem}

\appendix
\section{Connected components of $\FatP_\B$}

Here we turn to the proof of Lemma~\ref{lem:components}.  We will show
that in each connected component of $\FatP_\B$ there is, after making
a few controlled choices, a graph in a particular ``normal'' form, and
that this form is uniquely determined, up to these choices, by the
isomorphism type of its associated \ocobo.  This will show that the
connected components of $\FatP_\B$ are as desired.  Our proof of this
relies heavily on the corresponding construction for chord diagrams in
the closed string case, as presented by R.~Cohen and
V.~Godin~\cite{cohen-godin:polarized-view-string-topo}.

\begin{remark}
  We include this proof for completeness, but it will be superseded by
  a result stating that the category $\FatP_\B$ realizes to a space
  homotopy equivalent to an appropriate moduli space of \oc
  Riemann surfaces; this will recover Lemma~\ref{lem:components} upon
  applying $\pi_0$.
\end{remark}

\subsection{Isomorphism invariants of \osurfs}
%  Conversely to the first statement of
%  Proposition~\ref{prop:fattening-good}, we would like to show that two
%  well-partitioned \ofgraphs with isomorphic corresponding \ocobos are
%  in the same connected component of $\FatP_\B$. As a preliminary step
%  towards this, we describe a complete set of invariants of \bcobos
%  together with corresponding invariants of fat $\B$-graphs.

There are a few invariants which together determine the isomorphism
type of a connected \osurf.  They are as follows:
\begin{itemize}
\item The genus $g$ of $S$.
\item The subset $\bdryS S$ of $\pi_0(\bdrySt S)$ consisting of closed
  string boundary components.
\item The subset $\bdryI S$ of $\pi_0(\bdrySt S)$ consisting of string
  boundary intervals.
\item The function $\xi:\bdryI S\to \B$ assigning to a string boundary
  interval the label of its final endpoint, where ``final'' is with
  respect to the orientation.
\item A permutation $\psi\in\sym(\bdryI S)$ taking a string boundary
  interval $c$ to the one following $c$ in the same boundary
  component, in the direction induced by the orientation.
\item The subset $\bdryW S$ of $\pi_0(\bdryF S)$ consisting of closed
  components (``windows'').
\item The $\B$-labeling $\beta:\bdryW S\to \B$ induced from
  $\beta:\bdryF S\to\B$
\end{itemize}
If $S$ is an \ocobo, we can further identify:
\begin{itemize}
\item The partitions $\bdryS S = \bdryS^- S\disjcup\bdryS^+ S$,
  $\bdryI S = \bdryI^- S\disjcup\bdryI^+ S$.
\end{itemize}

Denote by $X(S)$ the tuple
\begin{displaymath}X(S):=(g,\bdryS S,\bdryI S,\xi,\psi,\bdryW S,
  \beta).\end{displaymath} For an \ocobo, denote by $Y(S)$ the tuple 
\begin{displaymath}Y(S):=(g,\bdryS S,\bdryI S,\xi,\psi,\bdryW S,
  \beta,\bdryS^\pm S, \bdryI^\pm S).\end{displaymath}
\begin{definition}
  In the absence of an \osurf or cobordism $S$, we consider tuples
  \begin{displaymath}(g,\bdryS,\bdryI,\xi,\psi,\bdryW,\beta)\end{displaymath}
  where $g\ge0$ is an integer, $\bdryS$, $\bdryI$ and $\bdryW$ are
  arbitrary finite sets, $\xi:\bdryI\to\B$ and $\beta:\bdryW\to\B$ are
  arbitrary functions, and $\psi$ is a permutation of $\bdryI$.  We
  call these \emph{\oc data tuples}. A \emph{partitioned} \oc data
  tuple is one of the form $(X,\bdryS^\pm,\bdryI^\pm)$, where
  $X=(g,\bdryS,\bdryI,\ldots)$ is an \oc data tuple and the
  $\bdryS^\pm$, $\bdryI^\pm$ are partitionings of $\bdryS$, $\bdryI$.

  Two such tuples (partitioned or not) are \emph{isomorphic} if they
  have the same $g$ and there are bijections of the corresponding sets
  which preserve the $\B$-labelings, the permutation $\psi$, and the
  partitions if present.
\end{definition}

\begin{proposition}
  \begin{enumerate}
  \item Two \osurfs (resp. cobordisms) are isomorphic if and only if
    their data tuples (resp. partitioned data tuples) are isomorphic. 
  
  \item Given an arbitrary (resp. partitioned) \oc data tuple, there
    is an \osurf $S$ with $X(S)$ isomorphic to it (resp. an \ocobo
    with $Y(S)$ isomorphic to it).
  \end{enumerate}
\end{proposition}
\begin{proofsketch}
  This is mostly clear, so we only provide a sketch of the
  construction for the second statement, which should make the first
  statement obvious.

  Choose an ordinary oriented surface $S$ of genus $g$ having its
  boundary components in bijection with $\bdryS\disjcup
  \bdryI/\gens{\psi}\disjcup\bdryW$. Label each entire component
  corresponding to $w\in\bdryW$ by $\beta(w)$, making it part of the
  free boundary.  Given a cycle $c=(x_1,\ldots,x_k)$ of $\xi$, let
  $A_c\subseteq\bdryS$ be the corresponding boundary component.
  Choose an embedding $\set{x_1,\ldots,x_k}\cross[0,1]\inclto A_c$
  which is orientation-preserving such that the cyclic order of the
  $x_i$ induced from the orientation of $A_c$ coincides with the
  cyclic order given by $\xi$. Declare the image of this embedding to
  belong to the string boundary. Decree the components of
  $\overline{A_c\setminus(\set{x_1,\ldots,x_k}\cross[0,1])}$ to be in
  the free boundary, and label them according to the rule that the
  component coming after $\set{x}\cross[0,1]$ in the cyclic order
  carries the label $\xi(x)$. The boundary components corresponding to
  $\bdryS$ are decreed to be part of the string boundary.  We omit the
  easy verification that this yields an \osurf having $X(S)\isom X$.
\end{proofsketch}

\bigskip%
Now, given an \ofgraph $\G$, $X(S(\G))$ is entirely determined by
$\G$, and in fact we can write $X(\G)$ for an \oc data tuple obtained
directly from $\G$, as follows:
\begin{displaymath}
  X(\G):=(g(\G),\bdryS\G,\bdryI\G,\xi,\psi,\bdryW\G,\restr{\beta_\G}{\bdryW\G}),
\end{displaymath}
where
\begin{itemize}
\item $\bdryI\G, \bdryS\G\subseteq \pi_0(\bdrySt\G)$ are the sets of
  free boundary intervals and cycles, respectively.
\item $\bdryW\G\subseteq \pi_0(\bdryF\G)$ is the set of closed free
  boundary components.
\item $g(\G) := 1-\frac{\#\pi_0(\bdry\G)+\chi(\G)}{2}$.
\item For $c\in\bdryI\G$, $\xi(c)$ is the $\B$-label carried by the
  final endpoint of $c$
\item $\psi$ takes a string boundary interval $c$ to the one
  appearing immediately after it in the component of $\bdry\G$
  containing $c$.
\end{itemize}
If $\G$ has a partitioning, we can further define
$Y(\G)=(X(\G),\bdryS^\pm\G,\bdryI^\pm I)$.  The following is clear.
\begin{proposition}
  $Y(S(\G))\isom Y(\G)$.
\end{proposition}

% \begin{proposition}
%   $X(\G_1\glue\G_2)\isom X(S(\G_1)\glue S(\G_2))$.
% \end{proposition}

% \begin{corollary}
%   $S(\G_1)\glue S(\G_2)\isom S(\G_1\glue\G_2)$.
% \end{corollary}

\subsection{Preliminary reductions}
\label{sec:prelim-reductions}

Let $\G\in\FatP_\B$ and let $b\subseteq\bdryF\G$ be a linear component
of the free boundary.
\begin{definition}
  Denote by $b_+\subseteq\bdrySt\G$ (resp. $b_-$) the string boundary
  component that appears immediately after $b$ (resp., before $b$) in
  $\bdryF\G$ according to the orientation.  We can classify $b$ as
  belonging to one of four types:
  \begin{itemize}
  \item say that $b\in B^-(\G)$ if $b_-\subseteq\bdry^-\G$ and
    $b_+\subseteq\bdry^-\G$,
  \item say that $b\in B^+(\G)$ if $b_-\subseteq\bdry^+\G$ and
    $b_+\subseteq\bdry^+\G$,
  \item say that $b\in B^{-+}(\G)$ if $b_-\subseteq\bdry^-\G$ and
    $b_+\subseteq\bdry^+\G$,
  \item say that $b\in B^{+-}(\G)$ if $b_-\subseteq\bdry^+\G$ and
    $b_+\subseteq\bdry^-\G$.
  \end{itemize}
\end{definition}
We now identify a convenient class of graphs containing a
representative of each connected component.  
\begin{definition}
  Call $\G$ \emph{special} if the following conditions hold:
  \begin{enumerate}
  \item Each linear component of $\bdryF\G$ consists of a single
    vertex.
  \item Each cyclic component of $\bdryF\G$ has exactly one edge.
  \item The image loop in $\G$ of each cyclic free boundary component
    $b$ is attached to a trivalent vertex; this vertex is therefore
    incident with the two half-edges forming the loop and with a third
    edge; the loop is then a \emph{balloon} attached to $\G$ by this
    edge.
  \item Each $b\in B^{-+}(\G)\cup B^{+-}(\G)\cup B^+(\G)$ is attached
    to a univalent vertex of $\G$
  \item Each $b\in B^-(\G)$ is attached to a bivalent vertex of $\G$.
  \item There are no other bivalent or univalent vertices.
  \end{enumerate}
\end{definition}
A special \ofgraph $\G$ then has distinguished $\B$-labeled bivalent
vertices, leaves with $\B$-labeled endpoint, and $\B$-labeled
balloons.

\begin{proposition}
  Any $\G\in\FatP_\B$ is connected to a special graph by a sequence of
  morphisms in $\FatP_\B$.
\end{proposition}
\begin{proof}
  The first two and last conditions can be attained by collapsing the
  image in $\G$ of a maximal forest in $\bdryF\G$; after that the
  three other conditions can be attained by expanding suitable
  vertices into trees.  
\end{proof}
\begin{remark}
  We cannot have a $b\in B^-(\G)$ attached to a \emph{univalent}
  vertex---this vertex would in turn be the endpoint of an edge having
  both its orientations in $\bdry^-\G$, and the partitioning would not
  be admissible.
\end{remark}
\subsection{Normal forms}
\label{sec:norm-forms}

In view of previous section, it is enough to show that two
\emph{special} graphs $\G,\G'\in\FatP_\B$ with $Y(\G)=Y(\G')$ are in
the same connected component.

Our strategy for finding suitable normal forms will be to use the
algorithm by R.~Cohen and V.~Godin
in~\cite{cohen-godin:polarized-view-string-topo}, henceforth referred
to as ``algorithm~$V$.''  The idea is, roughly, to run $\G$ through
this algorithm, having it work on the underlying \fgraph $U(\G).$
However, it will not be quite this simple, since we have to be careful
for two reasons:
\begin{itemize}
\item The algorithm in \cite{cohen-godin:polarized-view-string-topo}
  assumes that the starting graph is a chord diagram. Since we use the
  laxer admissibility condition from
  Definition~\ref{def:admissible-partitioning}, we will work to
  achieve this form; see Lemma~\ref{lem:chord-diag} below.
\item \Ofgraphs may have incoming and outgoing string boundary
  \emph{intervals} in addition to cycles.
\item Whenever the algorithm expands an edge of $U(\G)$, there is at
  least one way to expand the corresponding edge of $\G$ while
  keeping $\G$ special; however, when the algorithm collapses an edge
  of $U(\G)$, the corresponding edge of $\G$ may not be collapsable,
  since it may result in joining two components of the image of
  $\bdryF\G$.
\end{itemize}

To describe the different cases that arise when finding the normal
forms, we introduce some terms.
\begin{definition}
  Say that $\G\in\FatP_\B$ is \emph{clean} if it is special and,
  additionally, $\bdryF\G=B^{-+}(\G)\cup B^{+-}(\G)$.

  Given a special $\G\in\FatP_\B$, let $w(\G)$ be the \ofgraph
  obtained from $\G$ by
  \begin{enumerate}
  \item removing $B^+(\G)\cup B^-(\G)$ from $\bdryF\G$, 
  \item removing all the cyclic components of $\bdryF\G$, \emph{as
      well as} their image loops in $\G$,
  \item removing from $\G$ any leaves or bivalent vertices created by
    the previous two steps.
  \end{enumerate}

  Notice that $w(\G)$ inherits a boundary partitioning from $\G$,
  since the only removed free boundary intervals lie between string
  boundary components on the same side of the partitioning of $\G$.

  We define a \emph{weak string boundary component} of $\G$ to be a
  string boundary component of $w(\G)$, and we let $\bdry^\pm_w\G :=
  \bdry^\pm w(\G).$ 
\end{definition}

The result we aim for is as follows.
\begin{lemma}\label{lem:normal-forms}
  Let $\G\in\aFP_\B$ be connected.
  \begin{itemize}
  \item \emph{Case 1.} Suppose that $\G$ has a weak incoming boundary
    cycle. Then, $\G$ is connected in $\FatP_\B$ to an \ofgraph of
    the following form:
    \begin{center}
      \picbox{mp/fbglarge.3001}.
    \end{center}

  \item \emph{Case 2.} Suppose that $\G$ has no weak incoming
    boundary cycles, but it has more than one topological boundary
    component.  Then, $\G$ is connected to a graph of the form:
    \begin{center}
      \picbox{mp/fbglarge.3002}.
    \end{center}

  \item \emph{Case 3.} Suppose that $\G$ has no weak incoming boundary
    cycles, and that $\G$ has exactly one topological boundary
    component.  Then, $\G$ is connected to a graph of the form:
    \begin{center}
      \picbox{mp/fbglarge.3003}.
    \end{center}
  \end{itemize}

  The symbols used in the pictures are as follows:
  \begin{itemize}
  \item A serrated portion of an edge (\!\!\!\picbox{mp/fbglarge.3009})
    stands for a (possibly empty) sequence of bivalent vertices
    carrying elements of $B^-(\G)$; the triangles point towards the
    topological boundary component containing them.

  \item The symbol \picbox{mp/fbgsquare.3004} represents a (possibly
    empty) sequence \picbox{mp/fbglarge.3007} in which where each
    element \picbox{mp/fbgsquare.3006} stands for a structure of the
    form
    \begin{center}
      \picbox{mp/fbglarge.3008}.
    \end{center}
    Here, $c$ stands for an element of $B^{+-}(\G)$, $d$ for one of
    $B^{-+}(\G)$, and each $a_i$ stands for an element of $B^+(\G)$
    (the meaning of the serrated edge is as before).

  \item The symbol \picbox{mp/fbgsquare.3005} stands for a (possibly
    empty) sequence of balloons of the form\\[-2em]
    \begin{center}
      \picbox{mp/fbglarge.3010}
    \end{center}
    (that is, loops in the image of $\bdryF\G$ labeled by elements
    $w_i\in\B$.)

  \end{itemize}
  Note that those boundary components that contain a serrated edge are
  weak incoming boundary components, and the ones that contain a
  \picbox{mp/fbgsquare.3004} are topological boundary components which
  contain part of the outgoing boundary.

  In the first two cases, we will use terminology partially adapted
  from \cite{cohen-godin:polarized-view-string-topo}, as follows:
  \begin{itemize}
  \item The component marked $c_0$ in the pictures will be called the
    \emph{outer} component (this is called the big incoming circle in
    \cite{cohen-godin:polarized-view-string-topo}, but it is an
    outgoing component in our case~$2$).
  \item The topological boundary components in the top- and
    bottom-right quadrants in cases~1 and~2 will be called
    \emph{simple outgoing cycles}.
  \item The topological boundary component in the top-right quadrant
    will be called the \emph{complicated outgoing cycle}.
  \item The weak incoming cycles obtained by going clockwise around the
    small circles on the lower left of case~$1$ will be called
    \emph{small incoming cycles}.
  \end{itemize}

  The uniqueness of these normal forms is as follows:
  \begin{itemize}
  \item In case~$1$ we may choose an arbitrary weak incoming cycle
    $c_0$ to be the outer component, and we may choose an arbitrary
    topological boundary component containing part of the outgoing
    boundary to be the complicated outgoing cycle. Moreover, we can
    specify arbitrarily the order in which the simple outgoing cycles
    and the small incoming cycles occur in the cyclic ordering around
    the central vertex.
  \item In case~$2$ we may choose an arbitrary topological boundary
    component $c_0$ (necessarily containing part of the outgoing
    boundary) to be the outer component, and one to be the complicated
    outgoing cycle. As in case~$1$, we can specify the order of the
    simple outgoing cycles arbitrarily.
  \item In all cases, we can choose an arbitrary linear ordering
    (consistent with the underlying cyclic order) for the $\B$-labels
    to appear in each \picbox{mp/fbgsquare.3004} or
    \picbox{mp/fbglarge.3009} cluster.
  \item In all cases, we can choose an arbitrary order for the labels
    $w_i$ of the balloons appearing in the \picbox{mp/fbgsquare.3005}.
  \end{itemize}

  In each case, the normal form is uniquely determined, after the
  corresponding choices have been made, by the combinatorial data
  carried by the invariant $Y(\G)$, and therefore uniquely determined
  by the isomorphism type of $S(\G)$.
\end{lemma}

\medskip
\begin{lemma}\label{lem:extra-admissible}
  If $\G\in\FatP_\B$ has nonempty $\bdry^-\G$, then it is connected by
  morphisms to a special graph $\G'$ for which every edge of $w(\G')$
  has \emph{exactly} one of its orientations in the incoming string
  boundary.
\end{lemma}
\begin{proof}
  Assume without loss of generality that $\G$ is connected and
  special. We aim to get rid of edges of $w(\G)$ having two outgoing
  orientations (since $\G$ is well-partitioned, there are no edges
  having two incoming orientations). We will do this inductively, by
  showing that we can reduce the number of such ``bad'' edges via
  morphisms that introduce only ``good'' edges.

  Suppose there is at least one bad edge. Since $\G$ is connected,
  there is a vertex not in the image of $\bdryF\G$, and incident to
  both a bad edge and to at least one edge taking part in the incoming
  boundary. The cyclic ordering at this vertex hence looks like:
  \begin{center}
    \picbox{mp/fbglarge.3011},
  \end{center}
  where the symbols $\ominus$, $\oplus$ denote incoming and outgoing
  boundary, respectively. We can modify $\G$ in two steps, as follows:
  \begin{center}
    \begin{tabular}{ccccc}
      \picbox{mp/fbglarge.3011}  & $\from$ & 
      \picbox{mp/fbglarge.3012}  & $\to$ & 
      \picbox{mp/fbglarge.3013}  
    \end{tabular}
  \end{center}
  The graph in the middle maps to the other two graphs by obvious
  morphisms in $\FatP_\B$; the one on the left is $\G$ and the one on
  the right has one less bad edge than $\G$. We leave it to the reader
  to verify that this works even if the bad edge is a loop.
\end{proof}

\begin{lemma}\label{lem:chord-diag}
  If $\G\in\FatP_\B$ has nonempty $\bdry^-\G$, then it is connected
  by morphisms in $\FatP_\B$ to a special $\G'\in\FatP_\B$ for which
  $\bdry^- w(\G')$ is embedded in $w(\G)$ in such a way that the
  complement graph is a forest.
\end{lemma}
\begin{proof}
  We may assume that $\G$ is special and satisfies the condition in
  the conclusion of Lemma~\ref{lem:extra-admissible}. So the set of
  edges of $\bdry^-w(\G)$ is already embedded in the set of edges of
  $w(\G$); it remains to make the vertices embedded too. For every
  vertex $v$ not in the image of $\bdryF w(\G)$, the angles around $v$
  must alternate between incoming and outgoing when traversed
  according to the cyclic ordering at $v$.  Because of this, the
  following transformation (illustrated in the case that $v$ has
  valence 6) may be used to replace the vertex $v$ by a tree,
  resulting in $\G'\in\FatP_\B$ mapping to $\G$ by a morphism:
  \begin{center}
    \begin{tabular}{ccc}
      $\G$&&$\G'$ \\[-2em]
      \picbox{mp/fbgsquare.2003}&$\from$&\picbox{mp/fbgsquare.2004}. 
    \end{tabular}
  \end{center}\hspace{0em}\\[-2em]
  It is clear that applying this transformation to every vertex
  achieves the desired condition.  
  % (We could use a similar transformation at vertices in the image of
  % $\bdryF(\G)$ to achieve $\bdry^-\G$ embedded in $\G$ with a forest
  % complement, but we will make no use of that.)
\end{proof}

\medskip\noindent With these preliminary result in place, we are ready to prove our main
lemma.
\begin{proofof}{Lemma~\ref{lem:normal-forms}}
  We may assume that $\G$ is special. Notice that we may disregard the
  balloons throughout, since any two special \ofgraphs that differ
  only on the location of balloons along the outgoing boundary are
  connected by a sequence of morphisms; we leave this as an exercise.
  Thus, we can let them move around the graph arbitrarily, and we can
  collect them at the end to form a single \picbox{mp/fbgsquare.3005}
  at the correct location.

  Assume first that $\G$ is clean, so that $\G = w(\G)$.  By
  Lemma~\ref{lem:chord-diag}, we may assume that $\bdry^-\G$ is
  embedded in $\G$ with a forest complement. Following
  \cite{cohen-godin:polarized-view-string-topo}, we will call the
  edges of the complement ``ghost edges.''

  We can apply much of Algorithm~$V$ to $\G$, by using an arbitrary
  (open or closed) incoming boundary component $c_0$ of $\G$ instead
  of the ``big incoming component,'' and we may treat the other
  incoming boundary intervals much of the way as if they were incoming
  boundary cycles. We leave it to the reader to verify that $\G$ can
  be transformed, following the first steps of Algorithm~$V$, into a
  form in which exactly one non-univalent vertex~$v_0$ is incident
  with more than one ghost edge, and in which every incoming boundary
  component other than $c_0$ is connected by exactly one ghost edge to
  $v_0$.

  The next step in Algorithm~$V$ is to push the ``small incoming
  cycles'' all the way to the ``right'' in the cyclic ordering at
  $v_0$. This relies on these boundary components being cycles; thus
  we are not able to do it with the incoming boundary intervals.
  However, after the previous step, we may collapse the unique ghost
  edge attached to each of the incoming boundary intervals (other than
  $c_0$). After doing this, we create a structure of
  type~\picbox{mp/fbgsquare.3006} for each incoming interval. Call the
  resulting graph State~$\star$.

  Every time a structure of type~\picbox{mp/fbgsquare.3006} is
  created, we will treat it thereafter as a single unit, so its
  constituent edges will not be collapsed. We will call edges that are
  \emph{not} in any of these structures ``active;'' we will call a
  vertex active if all its incident edges are active.

  Now, suppose that $c_0$ is an incoming boundary \emph{cycle}.  Then,
  we can in fact continue applying Algorithm~$V$ to the end, and we
  obtain, after relocating structures of
  type~\picbox{mp/fbgsquare.3006} by means of expansion/collapse pairs
  if necessary, precisely the normal form in case~1.

  Suppose that we are in case~2, so that $c_0$ must be a incoming
  boundary interval.  After State~$\star$ there are no incoming
  boundary cycles, and we have collapsed all the incoming weak
  boundary intervals but one, namely $c_0$.  To proceed, collapse
  $c_0$ except for its two endpoints (which carry labels in
  $B^{-+}(\G)$ and $B^{+-}(\G)$ respectively); this is possible since
  $c_0$ is still embedded in $\G$. This creates a structure of
  type~\picbox{mp/fbgsquare.3006}. After this, we collapse a maximal
  active subtree of $\G$ . The result is a special \ofgraph with all
  its active half-edges belonging to $\bdry^+\G$, having more than one
  topological boundary cycle and exactly one active vertex .

  Change notation, letting $c_0$ be any topological boundary cycle,
  which will act as our ``big cycle'' from now on.  Using
  Lemma~\ref{lem:extra-admissible}, we can replace $\G$ by one in
  which each edge has exactly one of its orientations in $c_0$ (to do
  this, we introduce a temporary boundary partitioning on $\G$ which
  makes $c_0$ incoming and the remaining components outgoing, and then
  apply the lemma; here it is essential that $\G$ has more than one
  topological boundary component). Then, as we did in
  Lemma~\ref{lem:extra-admissible} for the incoming boundary, we can
  further replace $\G$ by a graph in which $c_0\to\G$ is an embedding
  and the complement of $c_0$ in $\G$ is a forest. From this point on,
  we apply Algorithm~$V$, treating $c$ as if it were the ``big
  incoming cycle,'' resulting in the required normal form.

  In case~3, we have after State~$\star$ and after collapsing the
  remaining incoming interval and a maximal active tree that $\G$ has
  one active vertex, a single topological component and all its active
  half-edges in $\bdry^+\G$. Since $\bdry\G$ is connected, we can
  ensure that all the labeled univalent vertices are contiguous by
  using a sequence of expansions/collapse pairs. The normal form for
  this case then follows directly from
  Lemma~\ref{lem:one-component-fg} below on the structure of \fgraphs
  with a single vertex and a single boundary component. (This case
  does not arise in \cite{cohen-godin:polarized-view-string-topo}
  because they do not consider \fgraphs with empty incoming boundary.)

  If $\G$ is not clean, then we apply the preceding procedure to
  $w(\G)$. Every time our algorithm expands a vertex of $w(\G)$ into
  an edge, the corresponding operation may be carried out on $\G$.
  When the algorithm contracts an edge of $w(\G)$, we must be careful
  because this edge might come from a sequence of edges in $\G$
  separated by bivalent vertices carrying labels in $B^-(\G)$.
  However, we can still collapse the corresponding edges in $\G$
  provided that we first move those bivalent vertices further along in
  the boundary of $\G$; this can be done by a straightforward sequence
  of expansion/collapse pairs. At the end, we are in a state for which
  $w(\G)$ is in one of the special forms; we can then use
  expansion/collapse pairs so that all the $B^-(\G)$ labels are
  contiguous (forming serrated edges), and so that all the univalent
  vertices carrying $B^+(\G)$ labels are also contiguous, and form
  part of a structure of type \picbox{mp/fbgsquare.3006}.

  We omit the proof of the uniqueness statement, which follows from a
  computation of the \oc data tuple of each of the normal form
  and from observing that they cover distinct isomorphism types.
\end{proofof}

\begin{lemma}\label{lem:one-component-fg}
  Let $\G$ be a \fgraph having a single boundary component and a
  single vertex. Then, $\G$ it is connected by a sequence of morphisms
  to a \fgraph having single vertex $v$ in which the cyclic ordering
  of half-edges at $v$ is of the form
  \begin{displaymath} (e_1,e_2,r(e_1),r(e_2),e_3,e_4,r(e_3),r(e_4),\ldots,e_{2k-1},e_{2k},r(e_{2k-1}),r(e_{2k}))  \end{displaymath}
  (where as usual $r$ stands for the edge-reversal involution on $\G$).
\end{lemma}
The proof is an easy induction.

\bigskip Lemma~\ref{lem:components} is now direct corollary of
Lemma~\ref{lem:normal-forms}.

\bibliographystyle{halpha}

%\bibliography{antonio}

\begin{thebibliography}{Ram05}

\bibitem[CG04]{cohen-godin:polarized-view-string-topo}
Ralph~L. Cohen and V{\'e}ronique Godin.
\newblock A polarized view of string topology.
\newblock In {\em Topology, geometry and quantum field theory}, volume 308 of
  {\em London Math. Soc. Lecture Note Ser.}, pages 127--154. Cambridge Univ.
  Press, Cambridge, 2004.

\bibitem[CJ02]{cohen:htpy-theoretic-realization}
Ralph~L. Cohen and John D.~S. Jones.
\newblock A homotopy theoretic realization of string topology.
\newblock {\em Math. Ann.}, 324(4):773--798, 2002.

\bibitem[CK]{cohen:generalized-collapse}
Ralph~L. Cohen and John~R. Klein.
\newblock To appear.

\bibitem[Cos]{costello:cy-cats}
Kevin~J. Costello.
\newblock {Topological conformal field theories and Calabi-Yau categories},
  arXiv:math.QA/0412149.

\bibitem[CS]{chas-sullivan:string-topology}
Moira Chas and Dennis Sullivan.
\newblock String topology.
\newblock {\em Annals of Math.}, arXiv:math.GT/9911159.
\newblock To appear.

\bibitem[Hara]{harrelson:thesis}
Eric Harrelson.
\newblock PhD thesis, University of Minnesotta.
\newblock In preparation.

\bibitem[Harb]{harrelson:homology-openclosed}
Eric Harrelson.
\newblock {On the homology of open-closed string field theory},
  arXiv:math.AT/0412249.

\bibitem[Har86]{harer:virtual-coho-dim-mcg}
John~L. Harer.
\newblock The virtual cohomological dimension of the mapping class group of an
  orientable surface.
\newblock {\em Invent. Math.}, 84(1):157--176, 1986.

\bibitem[LP]{lauda-pfeiffer:openclosed}
Aaron~D. Lauda and Hendryk Pfeiffer.
\newblock {Open-closed strings: Two-dimensional extended TQFTs and Frobenius
  algebras}.
\newblock AEI-2005-153; DAMTP-2005-80, arXiv:math.AT/0510664.

\bibitem[May99]{may:concise}
J.~P. May.
\newblock {\em A concise course in algebraic topology}.
\newblock Chicago Lectures in Mathematics. University of Chicago Press,
  Chicago, IL, 1999.

\bibitem[ML98]{maclane:categories}
Saunders Mac~Lane.
\newblock {\em Categories for the working mathematician}, volume~5 of {\em
  Graduate Texts in Mathematics}.
\newblock Springer-Verlag, New York, second edition, 1998.

\bibitem[Pen87]{penner:decorated-teichmuller-space}
R.~C. Penner.
\newblock The decorated {T}eichm\"uller space of punctured surfaces.
\newblock {\em Comm. Math. Phys.}, 113(2):299--339, 1987.

\bibitem[Ram05]{ramirez:thesis}
Antonio Ram\'{\i}rez.
\newblock {\em Open-Closed String Topology}.
\newblock PhD thesis, Stanford University, 2005.

\bibitem[Str84]{strebel:quadratic-differentials}
Kurt Strebel.
\newblock {\em Quadratic differentials}, volume~5 of {\em Ergebnisse der
  Mathematik und ihrer Grenzgebiete (3) [Results in Mathematics and Related
  Areas (3)]}.
\newblock Springer-Verlag, Berlin, 1984.

\bibitem[Sul04]{sullivan:openclosed}
Dennis Sullivan.
\newblock Open and closed string field theory interpreted in classical
  algebraic topology.
\newblock In {\em Topology, geometry and quantum field theory}, volume 308 of
  {\em London Math. Soc. Lecture Note Ser.}, pages 344--357. Cambridge Univ.
  Press, Cambridge, 2004.

\bibitem[Vor]{voronov:univ-algebra-notes}
Alexander~A. Voronov.
\newblock Notes on universal algebra, arXiv:math.QA/0111009.

\end{thebibliography}

\end{document}